\documentclass[a4paper, 11pt]{article}

\usepackage[utf8]{inputenc}

\usepackage{amsmath,amssymb,amsfonts}
\usepackage{mathtools}
\usepackage{booktabs}

\usepackage[usenames,dvipsnames]{xcolor}
\usepackage{graphicx}
\usepackage[labelformat=simple]{subcaption}

\usepackage{blindtext}

\usepackage{siunitx}
\DeclareSIUnit\flops{FLOP\per\second}
\DeclareSIUnit\flop{FLOP}
\DeclareSIUnit\ups{Up\per\second}
\DeclareSIUnit\lups{LUP\per\second}
\DeclareSIUnit\lu{LU}
\DeclareSIUnit\byte{B}
\DeclareSIUnit\bit{b}
\DeclareSIUnit\doubles{DP}
\DeclareSIUnit\cycles{cy}

\usepackage{tikz}
\usetikzlibrary{calc}
\usepackage{calc,pgf,pgfplots}
\usepackage{tikzscale}

\usepackage{listings}

\usepackage{calc}
\usepackage{algorithm,algpseudocode}
\usepackage{setspace}
\let\oldReturn\Return
\renewcommand{\Return}{\State\oldReturn}

\renewcommand{\Comment}[2][.5\linewidth]{%
  \leavevmode\hfill\makebox[#1][l]{$\triangleright$~#2}}

\usepackage{prettyref}
\newrefformat{tab}{Table~\ref{#1}}
\newrefformat{fig}{Figure~\ref{#1}}
\newrefformat{subfig}{Figure~\ref{#1}}
\newrefformat{alg}{Algorithm~\ref{#1}}
\newrefformat{sec}{Section~\ref{#1}}

\usepackage{geometry}
\lstdefinestyle{common-style}{
	basicstyle=\footnotesize,       % the size of the fonts that are used for the code
	numbers=left,                   % where to put the line-numbers
	numberstyle=\footnotesize,      % the size of the fonts that are used for the line-numbers
	stepnumber=1,                   % the step between two line-numbers. If it is 1 each line will be numbered
	numbersep=5pt,                  % how far the line-numbers are from the code
	showspaces=false,               % show spaces adding particular underscores
	showstringspaces=false,         % underline spaces within strings
	showtabs=false,                 % show tabs within strings adding particular underscores
	frame=single,                   % adds a frame around the code
	tabsize=2,                      % sets default tabsize to 2 spaces
	captionpos=b,                   % sets the caption-position to bottom
	breaklines=false,                % sets automatic line breaking
	breakatwhitespace=false,        % sets if automatic breaks should only happen at whitespace
	keywordstyle={\color{blue}\textbf},         % keywords are blue, (and blue)
	commentstyle={\color{OliveGreen}},          % comments
	literate={\$}{{\$}}1,
	escapechar=\&
}

\lstdefinestyle{fortran-style}{
	morecomment=[l][\color{\colorpragmaacc}]{!\$acc}, % ACC directives are orange
	morecomment=[l][\color{\colorpragmaomp}]{!\$omp},          % OMP directives are red
}

\usepackage{authblk}
\graphicspath{{.}}

\makeatletter
\newcommand{\settitle}{\@maketitle}
\makeatother
\newcommand{\hdg}{HDG}

\newcommand{\hdgseriestitle}{Linearizing the hybridizable discontinuous Galerkin method}

\newcommand{\of}[1]{\!\left(#1\right)}

\newcommand{\intelcompilerversion}{v.\,2018}

\newcommand{\name}[1]{\textsc{#1}}

\newcommand{\legendre}{\name{Legendre}}
\newcommand{\dirichlet}{\name{Dirichlet}}
\newcommand{\neumann}{\name{Neumann}}

\newcommand{\poisson}{\name{Poisson}}
\newcommand{\jacobi}{\name{Jacobi}}

\newcommand{\gllnames}{\name{Gau{\ss}-Lobatto-Legendre}}
\newcommand{\gausslobatto}{\name{Gau{\ss}-Lobatto}}
\newcommand{\lagrange}{\name{Lagrange}}
\newcommand{\gauss}{\name{Gau{\ss}}}
\newcommand{\schur}{\name{Schur}}
\newcommand{\schwarz}{\name{Schwarz}}

\newcommand{\galerkin}{\name{Galer\-kin}}

\newcommand{\domain}{\Omega}

\newcommand{\coordinatesymb}{x}
\newcommand{\coordinatesone}{\coordinatesymb_1}
\newcommand{\coordinatestwo}{\coordinatesymb_2}
\newcommand{\coordinatesthree}{\coordinatesymb_3}

\newcommand{\coordinates}{\physvector{x}}

\newcommand{\standardcoordinatesymb}{\xi}
\newcommand{\standardcoordinates}{\vec{\standardcoordinatesymb}}

\newcommand{\helmholtzparameter}{\lambda}

\newcommand{\physvector}[1]{\vec{#1}}

\newcommand{\normalvector}{\physvector{n}}
\newcommand{\normalvectorplus}{\normalvector^{\ +}}
\newcommand{\normalvectorminus}{\normalvector^{\ -}}

\newcommand{\elementindex}{e}
\newcommand{\nelement}{n_{\mathrm{\elementindex}}}
\newcommand{\nfaces}{n_{\mathrm{F}}}

\newcommand{\elementdomain}{\domain_{\elementindex}}
\newcommand{\standarddomain}{\domain^{\mathrm{S}}}
\newcommand{\poly}{p}

\newcommand{\eqdot}{\quad .}
\newcommand{\eqcomma}{\quad ,}
\newcommand{\varu}{u}
\newcommand{\varz}{z}
\newcommand{\varv}{v}
\newcommand{\varw}{\vec{w}}
\newcommand{\varf}{f}
\newcommand{\varmu}{\mu}

\newcommand{\varqsymb}{q}
\newcommand{\varq}{\varqsymb}
\newcommand{\varqvec}{\vec{\varqsymb}}
\newcommand{\varqvecplus}{\varqvec^{\ +}}
\newcommand{\varqvecminus}{\varqvec^{\ -}}

\newcommand{\mymat}[1]{\mathbf{#1}}

\newcommand{\vecu}{\mymat{u}}

\newcommand{\vecf}{\mymat{f}}
\newcommand{\vecF}{\mymat{F}}

\newcommand{\ndof}{n_{\mathrm{DOF}}}

\newcommand{\fluxu}{\tilde{\varu}}
\newcommand{\fluxq}{\vec{\tilde{\varqsymb}}}

\newcommand{\integd}{\,\mathrm{d}}

\newcommand{\domainintegral}[1]{%
  \int \limits _{\coordinates \in \domain}
  #1
  \integd \coordinates
}
\newcommand{\elementintegral}[1]{%
  \int \limits _{\coordinates \in \elementdomain}
  #1
  \integd \coordinates
}
\newcommand{\standardelementintegral}[1]{%
  \int \limits _{\standardcoordinatesymb \in \standarddomain}
  #1
  \integd \standardcoordinatesymb
}

\newcommand{\standardelementboundaryintegraloned}[1]{%
  \left[\left(#1\right)\of{-1} + \left(#1\right)\of{+1}\right]
}

\newcommand{\boundaryintegral}[1]{%
  \int \limits _{\coordinates \in \partial\domain}
  #1
  \integd \coordinates
}
\newcommand{\elementboundaryintegral}[1]{%
  \int \limits _{\coordinates \in \partial\elementdomain}
  #1
  \integd \coordinates
}
\newcommand{\neumannboundaryintegral}[1]{%
  \int \limits _{\coordinates \in \neumannboundary}
  #1
  \integd \coordinates
}

\newcommand{\neumannboundary}{\partial \domain_{\mathrm{N}}}
\newcommand{\dirichletboundary}{\partial \domain_{\mathrm{D}}}

\newcommand{\transmatsymb}{S}

\newcommand{\extransmat}{\mymat{\transmatsymb}}

\newcommand{\diffsymb}{D}
\newcommand{\diffmat}{\mymat{\diffsymb}}

\newcommand{\masssymb}{M}
\newcommand{\massmat}{\mymat{\masssymb}}

\newcommand{\stiffsymb}{L}
\newcommand{\stiffmat}{\mymat{\stiffsymb}}

\newcommand{\identitymat}{\mymat{I}}
\newcommand{\eigenvaluemat}{\mymat{\Lambda}}

\newcommand{\gathermatrix}{\mathbf{Q}^T}

\newcommand{\hdggathermatrix}{\mathbf{\tilde{Q}}^T}
\newcommand{\hdgscattermatrix}{\mathbf{\tilde{Q}}}

\newcommand{\tptwo}[2]{#1 \otimes #2}
\newcommand{\ptptwo}[2]{\left(#1 \otimes #2\right)}
\newcommand{\tp}[3]{#1 \otimes #2 \otimes #3}

\newcommand{\ptp}[3]{\left(#1 \otimes #2 \otimes #3\right)}
\newcommand{\geomcoeffsymb}{d}
\newcommand{\geomcoeff}[1]{\geomcoeffsymb_{#1}}

\newcommand{\diagonal}{\mymat{D}}

\newcommand{\npoint}{n_{\mathrm{p}}}
\newcommand{\order}[1]{\mathcal{O}\!\left(#1\right)}
\newcommand{\eigenspace}{\mathrm{E}}

\newcommand{\zeromat}{\mymat{0}}

\newcommand{\mata}{\mymat{A}}

\newcommand{\matb}{\mymat{B}}
\newcommand{\matc}{\mymat{C}}
\newcommand{\matd}{\mymat{D}}

\newcommand{\hdgk}{\mymat{K}}
\newcommand{\hdgkhat}{\mymat{\hat{K}}}

\newcommand{\hdgrhs}{\mymat{F}^{\mathrm{HDG}}}
\newcommand{\hdgrhshat}{\mymat{\hat{F}}^{\mathrm{HDG}}}
\newcommand{\hdgkelement}{\mymat{K}_e}
\newcommand{\hdga}{\mymat{A}}
\newcommand{\hdgaelement}{\mymat{A}_e}

\newcommand{\massmatstandard}{\massmat}
\newcommand{\hdguboundarystandard}{\hdge}
\newcommand{\hdguboundary}{\mymat{E}}

\newcommand{\hdgdiffmatstandard}{\hdgdiffmat}

\newcommand{\hdgfluxtouq}{\mymat{R}}
\newcommand{\hdguqtou}{\mymat{P}}

\newcommand{\hdgfluxtou}{\mymat{B}}
\newcommand{\hdgfluxtoustandard}{\hdgfluxtou}
\newcommand{\hdgfluxtoq}{\mymat{C}}
\newcommand{\hdgfluxtoqstandard}{\hdgfluxtoq}

\newcommand{\hdggstandard}{\hdgg}

\newcommand{\hdgutoq}{\mymat{D}}

\newcommand{\hdgsubmatrixqandflux}{\mymat{\tilde{A}}}

\newcommand{\hdgbs}{\mymat{B}_{\mathrm{S}}}
\newcommand{\hdgbssub}[1]{\mymat{B}_{\mathrm{S},#1}}

\newcommand{\hdgg}{\mymat{G}}

\newcommand{\hdggelement}{\hdgg_e}

\newcommand{\hdge}{\mymat{E}}

\newcommand{\hdgeelement}{\hdge_e}

\newcommand{\hdgzelement}{\mymat{Z}_e}
\newcommand{\hdgzelementthreed}{\mymat{Z}_{\mathrm{3D},e}}

\newcommand{\hdgtau}{\tau}
\newcommand{\hdgtauelement}{\hdgtau_e}
\newcommand{\hdgtauelementi}{\hdgtau_{e,i}}

\newcommand{\hdgtauplus}{\tau^{+}}
\newcommand{\hdgtauminus}{\tau^{-}}

\newcommand{\hdgtauhat}{\hat{\tau}}

\newcommand{\hdgdiffmat}{\diffmat}

\newcommand{\vecq}{\mymat{\varqsymb}}
\newcommand{\vecqvec}{\vec{\mymat{\varqsymb}}}
\newcommand{\vecutilde}{\mymat{\tilde{\varu}}}
\newcommand{\vecutildehat}{\mymat{\hat{\tilde{\varu}}}}
\newcommand{\vecRtilde}{\mymat{\tilde{r}}}
\newcommand{\vecRtildehat}{\mymat{\hat{\tilde{r}}}}

\newcommand{\vecztilde}{\mymat{\tilde{\varz}}}
\newcommand{\vecztildehat}{\mymat{\hat{\tilde{\varz}}}}

\newcommand{\solverhdgunprec}{HDG-unprec}
\newcommand{\solverhdgdiag}{HDG-diag}
\newcommand{\solverhdgblock}{HDG-block}
\newcommand{\solverhdgtrans}{HDG-trans}

\newcommand{\solverccgunprec}{CCG-unprec}
\newcommand{\solverccgdiag}{CCG-diag}
\newcommand{\solverccgblock}{CCG-block}
\newcommand{\solverccgtrans}{CCG-trans}

\author[a,b]{Immo Huismann\thanks{Corresponding author: Immo.Huismann@tu-dresden.de}}
\author[a,b]{Jörg Stiller}
\author[a,b]{Jochen Fröhlich}
\affil[a]{Institute of Fluid Mechanics, TU Dresden}
\affil[b]{Center for Advancing Electronics Dresden (cfaed)}
\title{\hdgseriestitle:~A linearly scaling operator}

\begin{document}
\maketitle
\begin{abstract}
  This paper proposes a matrix-free residual evaluation technique for the hybridizable discontinuous~\galerkin\ method requiring a number of operations scaling only linearly with the number of degrees of freedom.
  The method results from application of tensor-product bases on cuboidal Cartesian elements, a specific choice for the penalty parameter, and the fast diagonalization technique.
  In combination with a linearly scaling, face-wise preconditioner, a linearly scaling iteration time for a conjugate gradient method is attained.
  This allows for solutions in~\SI{1}{\micro\second} per unknown on one~CPU core -- a number typically associated with low-order methods.
\end{abstract}
\section{Introduction}\label{sec:intro}

At the front of spatial discretization, high-order methods are the current focus of research, especially continuous and discontinuous~\galerkin\ spectral-element methods~\cite{deville_2002_sem, karniadakis_1999_sem, hesthaven_2007_dg}.
Where low-order methods use linear shape functions or interpolations, these methods approximate the solution using piece-wise polynomial basis functions.
Using an ansatz of order~$\poly$, a convergence order of~${\poly+1}$ can be reached.
However, multiple issues arise:
First, while the number of degrees of freedom in three dimensions scales with~$\ndof \approx {\poly^3\nelement}$, where~$\nelement$ represents the number of elements, most operator evaluations scale with~${\order{\poly^4\nelement} = \order{\poly\ndof}}$ when exploiting tensor-product bases, i.e.\,super-linearly with the number of degrees of freedom~\cite{deville_2002_sem}.
% %
Second, the time step restriction from convection terms scales with~${1 / \poly^2}$ and even with~${1 / \poly^4}$ for diffusion~\cite{karniadakis_1999_sem}.
While fully explicit approaches are often preferred for compressible flows~\cite{hindenlang_2012_dg}, the incompressible case requires implicit treatment of diffusion and pressure to avoid crippling time stepping restrictions.
The occurring elliptic equations often take the form~${\lambda u - \Delta u = f}$ for a non-negative~$\lambda$~\cite{karniadakis_1991_projection, guermond_2006_projection}.
The runtime spent in the corresponding elliptic solvers can amount to~$90\ \%$ of the cost~\cite{fehn_2018_sim}.
Therefore, efficient solution of these equations is of central importance for a well-performing high-order flow solver.
The optimum would be a linearly scaling solver which, therefore, is the target.

This paper considers the discontinuous~\galerkin\ method~(DG) for discretization on Cartesian grids.
While~DG itself was first developed in conjunction with time-stepping for convection problems~\cite{jamet_1978_dg}, it was quickly adapted to spatial discretization~\cite{cockburn_1991_dg}, generating a plethora of different formulations ranging from interior penalty to local discontinuous ones~\cite{arnold_2002_dg}.
Here, the local discontinuous~\galerkin\ method~(LDG)~\cite{castillo_2000_ldg} is utilized.
The method allows for hybridization, i.e.\,using only the fluxes for the solution variable instead of the whole equation system.
This leads to the hybridized local discontinuous Galerkin method~(LDG-H), which reduces the number of algebraic unknowns from~$\order{\poly^3}$ to~$\order{\poly^2}$~\cite{kirby_2012_cg, yakovlev_2015_hdg, cockburn_2009_hdg}.
While every suboperator of the resulting equation system can be expressed in matrix-free form~\cite{kronbichler_2017_hdg}, the operator itself still scales with~${\order{\poly \ndof}}$.
For a solver to scale linearly with the number of unknowns when increasing the polynomial degree, a linearly scaling operator is mandatory.
While these linearly scaling operators are available for the continuous discretization~\cite{huismann_2017_condensation, huismann_2018_multigrid, huismann_2018_diss}, they have so far only been postulated for the discontinuous one~\cite{kronbichler_2017_hdg}, and, to the knowledge of the authors, no implementation thereof exists.

In previous work~\cite{huismann_2017_condensation, huismann_2018_multigrid}, the present authors have considered a continuous spectral-element discretization, developing a linearly scaling operator which, in conjunction with multigrid techniques, allows for a linearly scaling elliptic solver.
The resulting runtime scales with~${\order{\ndof}}$ independent of the polynomial degree.
The goal of this paper is to transfer these methods for the operator evaluation from the continuous discretization to the discontinuous one, i.e.\,to a hybridizable discontinuous~\galerkin\ operator.
It will be shown that, indeed, a linearly scaling operation account is achieved, which can be extended to a solver with an iteration time that scales linearly with the number of degrees of freedom and a runtime that scales linearly with them when increasing~$\poly$.
While the methods are presented here using~LDG-H, they can be similarly applied to the hybridizable interior penalty method described in~\cite{gander_2015_hdg}.

The structure of the paper is as follows:
First, the discretization with~LDG and~LDG-H is considered in~\prettyref{sec:operators}, and a linearly scaling operator derived for the one-dimensional case.
Then, the discretization and required operators are expanded to three dimensions in~\prettyref{sec:operators_threed}.
Thereafter in~\prettyref{sec:linearization}, sum and product factorization are applied to the operator, linearizing its operation count and validating the linear scaling with runtime tests.
Lastly, \prettyref{sec:runtimes} considers solvers with linearly scaling iteration time.

\section{The Hybridized Discontinuous Galerkin method in one dimension}\label{sec:operators}

\subsection{Local Discontinuous Galerkin discretization}

This paper considers the hybridized discontinuous~\galerkin\ method, introduced, for instance, in~\cite{kirby_2012_cg, yakovlev_2015_hdg, kronbichler_2017_hdg}, with a focus on the tensor-product operators.
To establish nomenclature and allow for a clearer understanding of the derivation of the linear solver, this section reiterates the basics thereof.
More expansive introductions to DG as well as~HDG can be found in~\cite{hesthaven_2007_dg, cockburn_2009_hdg, cockburn_2016_hdg}.

The considered elliptic equation reads
\begin{subequations}
  \begin{align}
\helmholtzparameter \varu - \Delta \varu &= \varf && \text{in}~\domain\label{eq:helmholtz_equation} \eqcomma
                                                                 \intertext{%
                                                              where~$\domain$ denotes the computational domain and, for each point~${\coordinates \in \domain}$, $\varu\of{\coordinates}$ the solution variable and~$\varf\of{\coordinates}$ the right-hand side.
                                                              Furthermore, $\helmholtzparameter$ constitutes a non-negative parameter.
                                                                 This equation is supplemented by boundary conditions
                                                                 }
                                                                 \varu &= g_{\mathrm{D}} && \text{on}~\dirichletboundary\\
    \normalvector \cdot \nabla\varu &= g_{\mathrm{N}} && \text{on}~\neumannboundary \eqcomma\label{eq:neumann_flux}
  \end{align}
\end{subequations}
with~$\dirichletboundary$ denoting a~\dirichlet\ boundary, $\neumannboundary$ a~\neumann\ boundary, and~$g_{\mathrm{N}}$ and~$g_{\mathrm{D}}$ the corresponding boundary values.
The introduction of an auxiliary variable~$\varqvec$ allows to rewrite~\eqref{eq:helmholtz_equation} in the so-called flux form
\begin{subequations}
  \begin{align}
    \helmholtzparameter \varu -\nabla \cdot \varqvec &= \varf \label{eq:helmholtz_flux_u}&&\\
    \varqvec &= \nabla \varu  \label{eq:helmholtz_flux_q} \eqcomma &&
               \intertext{with boundary conditions}
               \varu &= g_{\mathrm{D}} && \text{on}~\dirichletboundary\\
    \normalvector \cdot \varqvec &= g_{\mathrm{N}} && \text{on}~\neumannboundary \eqdot \label{eq:flux_neumann_flux}
  \end{align}\label{eq:helmholtz_flux}%
\end{subequations}
Using the test functions~$\varv, \varw, \varmu$, equations~\eqref{eq:helmholtz_flux_u}, \eqref{eq:helmholtz_flux_q}, and~\eqref{eq:flux_neumann_flux} can be cast into the weak form
\begin{subequations}
  \begin{align}
    \domainintegral{\helmholtzparameter \varu \varv} + \domainintegral{\nabla v \cdot \varqvec} &= \domainintegral{\varv\varf} + \boundaryintegral{\varv\normalvector \cdot \fluxq} \label{eq:helmholtz_weighted_form_with_fluxes_u} \\
    \domainintegral{\varw\cdot\varqvec} &= - \domainintegral{\nabla \cdot\varw \varu} + \boundaryintegral{\varw\cdot \normalvector \fluxu} \label{eq:helmholtz_weighted_form_with_fluxes_q}\\
    \boundaryintegral{\varmu\fluxq \cdot \normalvector} &= \neumannboundaryintegral{\varmu g_{\mathrm{N}}} \label{eq:flux_conservativity}   \eqcomma
  \end{align}\label{eq:helmholtz_weighted_form_with_fluxes}%
\end{subequations}
where the fluxes~$\fluxu$ and~$\fluxq$ approximate the boundary values of~$\varu$ and~$\varqvec$.

\begin{figure}[t]
  \hspace*{\fill}
  \begin{subfigure}{0.45\linewidth}
    \centering
    \includegraphics{hdg_element_variables.tikz}
    \caption{}\label{subfig:hdg_element_variables}
  \end{subfigure}
  \hfill
  \begin{subfigure}{0.45\linewidth}
    \centering
    \includegraphics{hdg_fluxes.tikz}
    \caption{}\label{subfig:hdg_fluxes}
  \end{subfigure}
  \hspace*{\fill}
  \caption{%
    Arrangement of variables for high order methods.
    Filled nodes indicate locations of the solution variables, whereas empty ones are only used in intermediate steps.
%    Left:~In the continuous case and primal~DG formulations the solution variable~$\varu$ suffices.
    %
    \subref{subfig:hdg_element_variables}:~In local discontinuous~DG formulations, the solution variable~$\varu$ and auxiliary variable~$\varqvec$ are utilized, with the flux only receiving an intermediate role.
    \subref{subfig:hdg_fluxes}:~In HDG, only the flux~$\vecutilde$ of the solution variable~$\varu$ remains.
  }%
  \label{fig:hdg_element_coupling}
\end{figure}
When discretizing the weighted flux form~\eqref{eq:helmholtz_weighted_form_with_fluxes_q}, the solution variables are~$\varu$ and~$\varqvec$, whereas the fluxes only serve as dependent variable used to couple the elements, as sketched in~\prettyref{fig:hdg_element_coupling}.
The main idea of~hybridization lies in eliminating~$\varu$, $\varqvec$, and~$\fluxq$ from system~\eqref{eq:helmholtz_weighted_form_with_fluxes} such that only the flux~$\fluxu$ remains as solution variable.
This leads to the condition that, after decomposing the domain~$\domain$ into~$\nelement$ elements~$\elementdomain$, the solution for~$\varu$ and~$\varqvec$ in each element should solely depend on the boundary flux~$\fluxu$ on each element, which limits the choice of the fluxes.
These are chosen in accordance with~\cite{cockburn_2016_hdg} here to be
\begin{subequations}
  \begin{align}
    \fluxu &= \frac{\hdgtauplus}{\hdgtauplus + \hdgtauminus} \varu^{+}
             + \frac{\hdgtauminus}{\hdgtauplus + \hdgtauminus} \varu^{-}
             - \frac{1}{\hdgtauplus + \hdgtauminus}\left(\varqvecplus \cdot \normalvectorplus + \varqvecminus \cdot \normalvectorminus\right) \label{eq:fluxu}\\
    \fluxq &= \frac{\hdgtauplus}{\hdgtauplus + \hdgtauminus} \varqvecplus
             + \frac{\hdgtauminus}{\hdgtauplus + \hdgtauminus} \varqvecminus
             - \frac{\hdgtauplus\hdgtauminus}{\hdgtauplus + \hdgtauminus}\left(\varu^{+} \cdot \normalvectorplus + \varu^{-} \cdot \normalvectorminus\right) \eqcomma
  \end{align}
\end{subequations}
where the superscripts~$+$ and~$-$ denote traces from the two adjoining elements on a boundary and~$\hdgtauplus$ and~$\hdgtauminus$ the corresponding penalty parameters.
For each element, the flux~$\fluxq$ equates to
\begin{align}
  \fluxq\of{\coordinates} &= \varqvec\of{\coordinates} - \hdgtauelement \normalvector\of{\coordinates} \left(\varu\of{\coordinates} - \fluxu\of{\coordinates}\right) \label{eq:flux_q} && \forall \coordinates \in \partial\elementdomain\eqcomma
\end{align}
where~$\hdgtauelement$ is the penalty parameter.

Inserting and~\eqref{eq:flux_q} into~\eqref{eq:helmholtz_weighted_form_with_fluxes} leads to the element equations.
This process simplifies~\eqref{eq:helmholtz_weighted_form_with_fluxes_u} to
\begin{subequations}
  \begin{align}
    &\begin{aligned}
      \elementintegral{\helmholtzparameter \varu \varv} + \elementboundaryintegral{\hdgtauelement\varv \varu} - \elementintegral{v \nabla \cdot \varqvec} &= \elementintegral{\varv\varf}
      +  \elementboundaryintegral{\hdgtauelement\varv \fluxu} \eqcomma
    \end{aligned}\label{eq:helmholtz_weighted_form_with_fluxes_u_two}
                                                                                                                                                            \intertext{%
                                                                                                                                                            again dropping the dependence on~$\coordinates$ for sake of readability.
                                                                                                                                                            Here, the third term was integrated by parts a second time to gain a symmetric equation system.
        Moreover, \eqref{eq:helmholtz_weighted_form_with_fluxes_q} becomes
        }
        &\elementintegral{\varw \cdot \varqvec} = - \elementintegral{\nabla \cdot\varw \varu} + \elementboundaryintegral{\varw \cdot \normalvector \fluxu}
        \intertext{%
        and, lastly, inserting~\eqref{eq:flux_q} into~\eqref{eq:flux_conservativity} yields
        }%
    &\sum_e\left[
    \elementboundaryintegral{\mu \varqvec \cdot \normalvector}
    +     \elementboundaryintegral{\mu \fluxu}
    -     \elementboundaryintegral{\hdgtauelement \mu \varu}
    \right] = \neumannboundaryintegral{\mu g_{\mathrm{N}}} \label{eq:transmission_condition} \eqdot
  \end{align}\label{eq:ldg_discretization}
\end{subequations}

\subsection{HDG~operators in one dimension}

In this section, a~HDG operator is derived in one dimension.
Here, $\varqvec$ has only got one component, $\varq$.
A set of polynomial basis functions~$\{\varphi_i\of{\standardcoordinatesymb}\}_{i=0}^{\poly}$ is introduced for~$\varu$ and~$\varq$ on the standard element~$\standarddomain = [-1,1]$ and mapped to each element.
Moreover, in each element~$\elementdomain$ $\fluxu$ has two values, each of them being one on their respective boundary while vanishing on the other, such that~two basis functions~$\phi_i:~\partial\standarddomain \to \{0,1\}$, $i \in \{1,2\}$, describe the fluxes~$\fluxu$ on the left and right boundary with~${\phi_1(-1) = 1,\phi_1(1) = 0}$ and~${\phi_2(-1) = 0,\phi_1(1) = 1}$.
Furthermore, the penalty parameter occuring in~\eqref{eq:flux_q} and~\eqref{eq:helmholtz_weighted_form_with_fluxes_u_two} is chosen as
\begin{align}
  \hdgtauelement &= \frac{2}{h_{e}} \hdgtauhat \eqcomma \label{eq:penalty_parameter_oned}
\end{align}
where~$h_e$ denotes the width of~$\elementdomain$.
Inserting the basis functions and penalty parameter into~\eqref{eq:ldg_discretization} leads to the following standard matrices and element matrices
\begin{subequations}
  \begin{align}
 \masssymb_{ij} &= \standardelementintegral{\varphi_i \varphi_j} & \massmat_e &= \frac{h_e}{2} \massmatstandard \\
  \diffsymb_{ij} &= \standardelementintegral{\varphi_i \partial_\standardcoordinatesymb \varphi_j} & \hdgdiffmat_e &= \hdgdiffmat \\
  E_{ij} &= +\hdgtauhat \standardelementboundaryintegraloned{\varphi_i \varphi_j} & \hdguboundary_e &= \frac{2}{h_e} \hdguboundarystandard\\
  G_{ij} &=  +\hdgtauhat \standardelementboundaryintegraloned{\phi_i \phi_j}  & \hdgg_e &= \frac{2}{h_e}\hdggstandard \\
  B_{ij} &= - \hdgtauhat \standardelementboundaryintegraloned{\varphi_i \phi_j} & \hdgfluxtou_e &= \frac{2}{h_e}\hdgfluxtoustandard \\
  C_{ij} &= \phantom{-\hdgtauhat}\standardelementboundaryintegraloned
           {\varphi_i \phi_j \normalvector} & \hdgfluxtoq_e &= \hdgfluxtoqstandard \eqdot
                                                              \intertext{Here, $\massmatstandard$ denotes the standard element mass matrix, $\hdgdiffmat$ the standard element differentiation matrix, and~$\stiffmat$ the stiffness matrix.
                                                              Furthermore, the stiffness matrix results from}
 \stiffmat &= \hdguboundary + \hdgdiffmat \massmatstandard^{-1} \hdgdiffmat^T & \stiffmat_e &= \frac{2}{h_e}\stiffmat \eqdot
\end{align}\label{eq:element_matrices}%
\end{subequations}
It is to be noted that the stiffness matrix~$\stiffmat$ differs from the one in the continuous case: The differentiation matrices are transposed and a further penalty term is present.

With the element matrices, the equation system~\eqref{eq:ldg_discretization} takes the form
\begin{align}
  \sum\limits_{e}
  \gathermatrix_e
  \begin{pmatrix}
    \helmholtzparameter\massmat_e + \hdguboundary_e & -\hdgdiffmat_e & \hdgfluxtou_e \\
    -\hdgdiffmat_e^T & -\massmat_e & \hdgfluxtoq_e \\
    \hdgfluxtou_e^T       & \hdgfluxtoq_e^T &  \hdggelement
  \end{pmatrix}
                                              \begin{pmatrix}
                                                \vecu_e\\\vecq_e\\\vecutilde_e
                                              \end{pmatrix}
                                                    &=
                                                      \sum\limits_{e}
                                                      \gathermatrix_e
                                                      \begin{pmatrix}
                                                        \massmat_e \vecf_e\\
                                                        \zeromat\\
                                                        \mymat{g}_{\mathrm{N},e}
                                                      \end{pmatrix}\label{eq:ldg_element_operator} \eqcomma
\end{align}
where~$\gathermatrix_e$ gathers the contributions from several elements for the global degrees of freedom and boldface denotes the coefficient vectors, e.g.\,$\vecu_e$ for the coefficients on the element for~$\varu$.
As the solution and the auxiliary variable are discontinuous, the element coefficients are the global degrees of freedom.
Therefore, the first two lines remain local to the element.
The last line couples the elements, with the right-hand side~$\mymat{g}_{\mathrm{N},e}$ denoting either the~\neumann\ boundary condition, or zero.
Eliminating~$\vecu_e$ and~$\vecq_e$ from~\eqref{eq:ldg_element_operator} leads to a global equation system of the form
\begin{align}
  \hdgk \vecutilde &= \hdgrhs
                     \intertext{which can be rewritten as}
  \sum_e \hdggathermatrix_e \hdgkelement \underbrace{\hdgscattermatrix_e  \vecutilde}_{\vecutilde_e} &= \sum_e \hdggathermatrix_e \hdgrhs_e \label{eq:hdg_system}
\end{align}
where~$\hdgkelement$ denotes the element operator, $\hdgrhs$ the right-hand side and~$\hdgscattermatrix_e$ the matrix mapping global flux degrees of freedom to element-local ones.
The element operators~$\hdgkelement$ result via~\schur\ complement of~\eqref{eq:ldg_element_operator}
\begin{align}
  \hdgkelement &= \hdggelement -
                 \underbrace{\begin{pmatrix}
                     \hdgfluxtou_e\\\hdgfluxtoq_e
                   \end{pmatrix}^T}_{\hdgfluxtouq_e^T}
  \hdgaelement^{-1}
  \underbrace{\begin{pmatrix}
      \hdgfluxtou_e\\\hdgfluxtoq_e
    \end{pmatrix}}_{\hdgfluxtouq_e}  \label{eq:hdg_operator_oned}\eqcomma
  \intertext{where~$\hdgaelement$ denotes the coupling between~$\vecu$ and~$\vecq$ such that}
    \hdgaelement &=
                   \begin{pmatrix}
                     \helmholtzparameter\massmat_e + \hdguboundary_e  & -\hdgdiffmat_e\\
                     -\hdgdiffmat_e^T & -\massmat_e
                   \end{pmatrix}\eqdot
                                        \intertext{The above matrix can be explicitely inverted, e.g. via~\textsc{Schur} complement, leading to}
                                        % \hdgaelement^{-1} &=
                                        %                     \begin{pmatrix}
                                        %                       \hdgzelement & -\hdgzelement \hdgdiffmat_e \massmat_e^{-1} \\
                                        %                       -{\massmat_e^{-T}\hdgdiffmat_e^{-T} \hdgzelement^{T}} &  - \massmat_e^{-1} \left(\identitymat - \hdgdiffmat_e^{T} \hdgzelement \hdgdiffmat_e \massmat^{-1}_e\right)
                                        %                     \end{pmatrix} \\
    \hdgaelement^{-1} &=
                      \underbrace{\begin{pmatrix}
                          \identitymat \\
                          - \massmat_e^{-T}\hdgdiffmat_e^{T}
                        \end{pmatrix}}_{\hdguqtou^T_e}
  \underbrace{\left(\lambda \massmat_e + \hdgeelement + \hdgdiffmat_e \massmat_e^{-1} \hdgdiffmat_e^T \right)^{-1}}_{\hdgzelement}
  \underbrace{\begin{pmatrix}
      \identitymat &
      - \hdgdiffmat_e\massmat_e^{-1}
    \end{pmatrix}}_{\hdguqtou_e}
                     -
                     \begin{pmatrix}
                       \zeromat & \zeromat\\
                       \zeromat & \massmat_e^{-1}
                     \end{pmatrix}\label{eq:element_operator_u_q_1d} \eqdot
\end{align}
The application of the element operator~$\hdgkelement$ contains four steps:
First, the diagonal matrix~$\hdggelement$ integrates on the faces.
Then, $\hdguqtou_e$ computes the residual induced by the fluxes~$\vecutilde_e$ in the inner element for~$\vecu_e$ and~$\vecq_e$, and from these the resulting~$\vecu_e$ and~$\vecq_e$ via~$\hdgaelement^{-1}$.
Lastly, applying~$\hdguqtou_e^T$ computes the effects of the solution infered into the element onto the flux~$\vecutilde_e$.

\subsection{Linearly scaling, matrix-free HDG operator evaluation}

For evaluation of the~HDG element operator, fast application of~$\hdgaelement^{-1}$ is key.
The matrix~$\hdgzelement$ is the main obstacle, as it is a dense, of size~${(\poly+1)\times (\poly+1)}$, and differs from element to element.
The specific choice of~$\hdgtauelement$ in~\eqref{eq:penalty_parameter_oned} allows to describe~$\hdgzelement$ with the same two matrices in every element:
\begin{align}
  \hdgzelement^{-1}&=\lambda \massmat_e + \hdgeelement + \hdgdiffmat_e \massmat_e^{-1} \hdgdiffmat_e^T\\
  \stiffmat_e &= \hdguboundary_e + \hdgdiffmat_e \massmat_e^{-1} \hdgdiffmat_e^T \label{eq:hdg_pseudo_stiffness_mat}\\
  \Rightarrow \hdgzelement   &= {\left( \lambda\massmat_e + \stiffmat_e \right)}^{-1} =  {\left( \frac{\lambda h_e}{2}  \massmatstandard + \frac{2}{h_e}\stiffmat \right)}^{-1} \eqdot
\end{align}
This allows for a generalized eigenvalue decomposition
\begin{subequations}
  \begin{align}
    \extransmat^T \massmatstandard \extransmat &= \identitymat \label{eq:generalized_eigenvalue_decomposition_massmat}\\
    \extransmat^T \stiffmat   \extransmat &= \eigenvaluemat \eqcomma
  \end{align}\label{eq:generalized_eigenvalue_decomposition}%
\end{subequations}%
where~$\extransmat$ is a non-orthogonal transformation matrix, and~$\eigenvaluemat$ a diagonal matrix containing the eigenvalues.
Using these, the inverse computes to
\begin{align}
  \hdgzelement &= \extransmat \diagonal_{\hdgzelement}^{-1} \extransmat^T \eqcomma
                 \intertext{where}
                 \diagonal_{\hdgzelement} &=  \frac{\lambda h_e}{2} \identitymat + \frac{2}{h_e}\eigenvaluemat \eqdot
\end{align}
With this representation, the inverse of~$\hdgaelement$ can be modified with respect to~\eqref{eq:element_operator_u_q_1d}, now reading
\begin{align}
  \hdgaelement^{-1} &=
  \hdguqtou^T_e\extransmat \diagonal_{\hdgzelement}^{-1} \extransmat^T\hdguqtou_e
                     -
                     \begin{pmatrix}
                       \zeromat & \zeromat\\
                       \zeromat & \massmat_e^{-1}
                     \end{pmatrix}\label{eq:element_operator_u_q_1d_with_s} \eqdot
\end{align}
Inserting~\eqref{eq:element_operator_u_q_1d_with_s} into the element operator~$\hdgkelement$ leads to
\begin{align}
    \hdgkelement &=
       \hdggelement
                   -
                   \underbrace{\begin{pmatrix}
                     \hdgfluxtou_e^T&
                     \hdgfluxtoq_e^T
                   \end{pmatrix}}_{\hdgfluxtouq_e^T}
    \left[
                     \hdguqtou_e^T
      \extransmat \diagonal_{\hdgzelement}^{-1} \extransmat^T
                     \hdguqtou_e
      -
      \begin{pmatrix}
        \zeromat & \zeromat\\
        \zeromat & \massmat_e^{-1}
      \end{pmatrix}
    \right]
                   \underbrace{\begin{pmatrix}
                     \hdgfluxtou_e \\ \hdgfluxtoq_e
                   \end{pmatrix}}_{\hdgfluxtouq_e}
  \eqcomma
   \intertext{such that further evalation yields}
    \hdgkelement &= \hdggelement
                   + \hdgfluxtoq_e^T  \massmat_e^{-1}\hdgfluxtoq_e
                   -
                   \hdgfluxtouq_e^T\hdguqtou_e^T\extransmat
                   \diagonal_{\hdgzelement}^{-1}\
                   \extransmat^T\hdguqtou_e\hdgfluxtouq_e\\
  \intertext{or, equivalently,}
    \hdgkelement
                   &= \frac{2}{h_e} \hdggstandard
                   + \frac{2}{h_e} \hdgfluxtoq^T  \massmat^{-1}\hdgfluxtoq
                   - \frac{2}{h_e}\hdgfluxtouq^T\hdguqtou^T\extransmat
                   \diagonal_{\hdgzelement}^{-1}\
                   \extransmat^T\hdguqtou_{\phantom{e}}\hdgfluxtouq_{\phantom{e}} \frac{2}{h_e} \label{eq:hdgkelement_factorized} \eqdot
  \end{align}

  The last form consists of three terms.
  Application of the first two just requires one $2\times2$ matrix product, whereas the latter consists of using a row matrix to generate the right-hand side for~$\vecu$, using the inverse eigenvalues, and then a reduction back to the fluxes.

  The ultimate goal lies in a linearly scaling, matrix-free~HDG operator in three dimensions.
  To this end, the one-dimensional case of~\eqref{eq:hdgkelement_factorized} serves as a prototype.
  To attain a linearly scaling algorithm for the residual evaluation, the row matrix mapping into the element eigenspace requires an explicit representation:
  \begin{align}
     \hdgbs &= \extransmat^T\hdguqtou \hdgfluxtouq  = \extransmat^T
                                             \begin{pmatrix}
                                               \identitymat & - \hdgdiffmat\massmat^{-1}
                                             \end{pmatrix}
                                                              \begin{pmatrix}
                                                                \hdgfluxtou \\
                                                                \hdgfluxtoq
                                                              \end{pmatrix}
      =
      \extransmat^T\hdgfluxtou
      -
    \extransmat^T\hdgdiffmatstandard \massmat^{-1} \hdgfluxtoq \label{eq:hdgbs} \eqdot
    \end{align}
    \begin{algorithm}[t]
      \caption{%
        Computation of the effect of the HDG operator on a per-element basis, where~${\geomcoeff{1,e} = h_e/2}$ denotes the metric coefficient.
      }%
      \label{alg:hdg_op_1d}
      \begin{algorithmic}
        \Function{HDG\_Op\_1D}{$\vecutilde$}
        \For{$ e = 1, \nelement$}
        \State{$\vecF_{\eigenspace} \gets \geomcoeff{1,e}\hdgbs\vecutilde_{e}$}
        \State{$\vecu_{\eigenspace} \gets \diagonal_{\hdgzelement}^{-1} \vecF_{\eigenspace}$}
        \State{$\vecRtilde_{e} \gets \geomcoeff{1,e}(\hdgg  + \hdgfluxtoq^T \massmatstandard^{-1}\hdgfluxtoq ) \vecutilde_{e} - \geomcoeff{1,e}\hdgbs ^T \vecu_{\eigenspace}$}
        \EndFor{}
        \Return $\vecRtilde$
        \EndFunction
      \end{algorithmic}
    \end{algorithm}%
    Using this mapping, the element-wise~\hdg\ residual can be evaluated via~\prettyref{alg:hdg_op_1d}, which requires~$\order{\poly\nelement}$ operations.
    First, the operator computes the residual~$\vecF_{\eigenspace}$ in the element eigenspace by~applying~$\hdgbs$ to the fluxes and then computes a solution~$\vecu_{\eigenspace}$ in the eigenspace by applying the inverse eigenvalues~$\diagonal_{\hdgzelement}^{-1}$.
    Lastly, the effect of~$\vecu_{\eigenspace}$ onto~$\vecutilde_e$ results by applying~$\hdgbs^T$ on the solution~$\vecu_{\eigenspace}$ and the diagonal term.
    The three operators of mapping into the element eigenspace, applying the inverse eigenvalues, and, lastly, mapping back generate an operation count~$5 (\poly+1)$, whereas the application of~${(\hdgg  + \hdgfluxtoq^T \massmatstandard^{-1}\hdgfluxtoq )}$ requires~$8$ operations per element.
    Therefore, \prettyref{alg:hdg_op_1d} constitutes a linearly scaling, matrix-free residual evaluation technique, albeit only for the one-dimensional case.

\section{The~HDG element operator in  three dimensions}\label{sec:operators_threed}

\subsection{Tensor-product matrices}

In this section, a tensor-product basis is considered.
This provides further structure to the operators and leads to tensor-product matrices, as introduced in~\cite{lynch_1964_tensors, deville_2002_sem}.
For~$\mata,\matb\in \mathbb{R}^{n,n}$, the tensor product~${\tptwo{\matb}{\mata} \in \mathbb{R}^{n^2,n^2}}$ is defined as
\begin{align}
  \tptwo{\matb}{\mata} &=
                         \begin{pmatrix}
                           \mata B_{11} & \mata B_{12} & \dots  &  \mata B_{1n} \\
                           \mata B_{21} & \mata B_{22} & \dots  &  \mata B_{2n} \\
                           \vdots       & \vdots       & \ddots & \vdots \\
                           \mata B_{n1} & \mata B_{n2} & \dots  & \mata B_{nn}
                         \end{pmatrix} \eqdot \label{eq:tp_definition}
\end{align}
Where the application of the matrix~$\tptwo{\matb}{\mata}$ as a whole to a vector~$\vecu \in \mathbb{R}^{n^2}$ incurs~${2n^4}$ operations, the tensor-product can be applied by, first, using matrix~$\mata$ in the first direction and then $\matb$ in the second.
As only two one-dimensional matrix products are involved, ${2 \cdot 2n^{2+1}}$ operations occur.
Furthermore, the following properties directly result from~\eqref{eq:tp_definition}:
\begin{subequations}
  \begin{align}
    {\ptptwo{\matb}{\mata}}^T &= \tptwo{\matb^T}{\mata^T}\\
    {\ptptwo{\matb}{\mata}}{\ptptwo{\matd}{\matc}} &= \tptwo{(\matb \matd)}{(\mata\matc)} \label{eq:tp_property_fusing} \eqdot
  \end{align}%
\end{subequations}
The extension to three dimensions is straight-forward via~$\tptwo{\matc}{(\tptwo{\matb}{\mata})}$, retaining the properties and raising the operation count to~${3\cdot 2 n^{3+1}}$, as is the extension to non-square matrices.

\subsection{LDG operator structure in three dimensions}
\prettyref{alg:hdg_op_1d} provides a linearly scaling residual evaluation for~HDG in one dimension.
In this section it will be extended to three dimensions by inferring a tensor-product structure into the three-dimensional~HDG operator.

When expanding the~LDG discretization to three dimensions, the vector~$\vecqvec$ has three components, $\vecq_{1}$, $\vecq_{2}$, and~$\vecq_{3}$.
Here, only Cartesian cuboidal elements with tensor-product bases are considered, which allows for a separation of the dimensions:
The flux vector~$\vecutilde$ can be separated into three parts, $\vecutilde_{1}$, $\vecutilde_{2}$, and~$\vecutilde_{3}$, each residing on the faces in the respective direction.
Moreover, the axis-aligned grid removes inter-dependencies in the derivatives:
In each direction~$\coordinatesymb_i$, the auxiliary variable~$\vecq_{i}$ solely depends on~$\vecu$ and the corresponding flux~$\vecutilde_{i}$.
This yields the element equation system
%~\eqref{eq:ldg_element_operator} directly extends to the three-dimensional case in each dimension
\begin{align}
  \sum\limits_{e}\gathermatrix_e
   \begin{pmatrix}
    \lambda \massmat_e +  \hdgeelement & - \hdgdiffmat_{1,e}& - \hdgdiffmat_{2,e}& - \hdgdiffmat_{3,e} & \hdgfluxtou_{1,e}       & \hdgfluxtou_{2,e}       & \hdgfluxtou_{3,e}      \\
     - \hdgdiffmat_{1,e}^T       & -\massmat_e & \zeromat   & \zeromat   & \hdgfluxtoq_{1,e} & \zeromat & \zeromat  \\
     - \hdgdiffmat_{2,e}^T       & \zeromat   & -\massmat_e & \zeromat   & \zeromat & \hdgfluxtoq_{2,e} & \zeromat  \\
     - \hdgdiffmat_{3,e}^T     & \zeromat   & \zeromat   & -\massmat_e & \zeromat & \zeromat & \hdgfluxtoq_{3,e} \\
     \hdgfluxtou_{1,e}  ^T  & \hdgfluxtoq_{1,e}^T & \zeromat & \zeromat   & \hdgg_{1,e} & \zeromat   & \zeromat   \\
     \hdgfluxtou_{2,e}  ^T  & \zeromat & \hdgfluxtoq_{2,e}^T & \zeromat   & \zeromat   & \hdgg_{2,e} & \zeromat   \\
     \hdgfluxtou_{3,e}^T  & \zeromat & \zeromat & \hdgfluxtoq_{3,e}^T & \zeromat   & \zeromat   & \hdgg_{3,e}
   \end{pmatrix}
   \begin{pmatrix}
     \vecu_e \\ \vecq_{1,e} \\ \vecq_{2,e} \\ \vecq_{3,e} \\ \vecutilde_{1,e} \\ \vecutilde_{2,e} \\ \vecutilde_{3,e}
   \end{pmatrix}
                                          &=
                                            \sum\limits_{e}\gathermatrix_e
    \begin{pmatrix}
      \massmat_e \vecf_e\\ \zeromat\\ \zeromat\\ \zeromat  \\ \mymat{g}_{1,\mathrm{N},e} \\ \mymat{g}_{2,\mathrm{N},e}  \\  \mymat{g}_{3,\mathrm{N},e}
    \end{pmatrix}  \eqcomma                                          \label{eq:ldg_element_operator_threed}
\end{align}
where the subscript~$e$ denotes the element operator for the three-dimensional element, and the subscripts~$1$, $2$, and~$3$ indicate the respective direction for the operation.
For application of the matrix in~\eqref{eq:ldg_element_operator_threed}, similar operators as in the one-dimensional case are required.
Let~$h_{i,e}$ denote the element width of~$\elementdomain$ in direction~$\coordinatesymb_i$, which leads to the metric factors
\begin{align}
  \begin{pmatrix}
    \geomcoeff{0,e}&  \geomcoeff{1,e}&\geomcoeff{2,e}&\geomcoeff{3,e}
  \end{pmatrix}
  &=
    \frac{h_{1,e}h_{2,e}h_{3,e}}{8}
    \begin{pmatrix}
      1 & {\left(\frac{2}{h_{1,e}}\right)}^2 & {\left(\frac{2}{h_{2,e}}\right)}^2 & {\left(\frac{2}{h_{3,e}}\right)}^2
    \end{pmatrix} \eqdot
                                                                                    \intertext{As in the one-dimensional case, the penalty parameter is chosen as}
    \hdgtau_{i,e} &= \frac{2}{h_{i,e}} \hdgtauhat \eqcomma
\end{align}
which allows for factorization further down the line.
A three-dimensional tensor-product base~${\{\varphi_i\of{\standardcoordinatesymb_1}\varphi_j\of{\standardcoordinatesymb_2}\varphi_k\of{\standardcoordinatesymb_3}\}}_{i,j,k=0}^{\poly}$ on the standard element~${\standarddomain = {[-1,1]}^3}$ is utilized, where~$\standardcoordinates$ denotes the standard coordinates in three dimensions.
Furthermore, the basis for the fluxes results from the tensor-product base at the element boundaries, e.g.\,$\standardcoordinatesymb_1 = -1$.
\prettyref{fig:hdg_element_coupling} depicts the two-dimensional case.
Using this notation, the insertion of the tensor-product base into~\eqref{eq:ldg_discretization} leads to the tensor-product matrices listed in~\prettyref{tab:hdg_matrices_threed}.

\begin{table}[t]
  \caption{%
    System matrices occurring in a three-dimensional~LDG formulation on Cartesian tensor-product elements expressed using the standard element matrices stemming from~\eqref{eq:element_matrices}.
    For the tensor-product structure, the vectors~$\vecutilde_{i,  e}$ have the shape~${2\times(\poly+1)\times(\poly+1)}$, ${(\poly+1)\times2\times(\poly+1)}$, and~${(\poly+1)\times(\poly+1)\times2}$ for the~$\coordinatesone$, $\coordinatestwo$, and~$\coordinatesthree$ direction, respectively.
  }
  \centering
    \begin{tabular}{llrll|lllr}\toprule
      Matrix         && Definition &&&& Matrix         && Definition\\\midrule
      $\massmat_e $          && $\geomcoeff{0,e}\tp{\massmatstandard}{\massmatstandard}{\massmatstandard}$                    &&&&$\hdgeelement$         && $\hdguboundary_{1,e} + \hdguboundary_{2,e} + \hdguboundary_{3,e}$\\
      $\hdgdiffmat_{1,e} $   && $h_{1,e}/2\ \geomcoeff{1,e}\tp{\massmatstandard}{\massmatstandard}{\hdgdiffmatstandard}  $ &&&& $\hdguboundary_{1,e}$  && $\geomcoeff{1,e} \tp{\massmatstandard}{\massmatstandard}{\hdguboundarystandard}$\\
      $\hdgdiffmat_{2,e} $   && $h_{2,e}/2\ \geomcoeff{2,e}\tp{\massmatstandard}{\hdgdiffmatstandard}{\massmatstandard}  $ &&&& $\hdguboundary_{2,e}$  && $\geomcoeff{2,e} \tp{\massmatstandard}{\hdguboundarystandard}{\massmatstandard}$\\
      $\hdgdiffmat_{3,e} $   && $h_{3,e}/2\ \geomcoeff{3,e}\tp{\hdgdiffmatstandard}{\massmatstandard}{\massmatstandard}  $ &&&& $\hdguboundary_{3,e}$  && $\geomcoeff{3,e} \tp{\hdguboundarystandard}{\massmatstandard}{\massmatstandard}$\\\midrule

      $\hdgg_{1, e}$         && $\geomcoeff{1,e} \tp{\massmatstandard}{\massmatstandard}{\hdggstandard}$          &&&& \\
      $\hdgg_{2, e}$         && $\geomcoeff{2,e} \tp{\massmatstandard}{\hdggstandard}{\massmatstandard}$          &&&& \\
      $\hdgg_{3, e}$         && $\geomcoeff{3,e} \tp{\hdggstandard}{\massmatstandard}{\massmatstandard}$          &&&& \\
      $\hdgfluxtou_{1,e}$  && $\geomcoeff{1,e}\tp{\massmatstandard}{\massmatstandard}{{\hdgfluxtoustandard}}$ &&&& $\hdgfluxtoq_{1,e}$  && $h_{1,e}/2\ \geomcoeff{1,e}  \tp{\massmatstandard}{\massmatstandard}{{\hdgfluxtoqstandard}} $ \\
      $\hdgfluxtou_{2,e}$  && $\geomcoeff{2,e}\tp{\massmatstandard}{{\hdgfluxtoustandard}}{\massmatstandard}$ &&&& $\hdgfluxtoq_{2,e}$  && $h_{2,e}/2\ \geomcoeff{2,e}  \tp{\massmatstandard}{{\hdgfluxtoqstandard}}{\massmatstandard} $ \\
      $\hdgfluxtou_{3,e}$  && $\geomcoeff{3,e}\tp{{\hdgfluxtoustandard}}{\massmatstandard}{\massmatstandard}$ &&&& $\hdgfluxtoq_{3,e}$  && $h_{3,e}/2\ \geomcoeff{3,e}  \tp{{\hdgfluxtoqstandard}}{\massmatstandard}{\massmatstandard} $ \\\bottomrule
    \end{tabular}%
  \label{tab:hdg_matrices_threed}
\end{table}

\subsection{HDG~operator in three dimensions}

Equation~\eqref{eq:ldg_element_operator_threed} defines the element contributions to the global~LDG operator.
Application of a~\schur\ complement eliminates~$\vecu$ and~$\vecqvec$ on a per-element basis and leads to an equation system of the form~\eqref{eq:hdg_system} with the element operator becoming
\begin{align}
  \hdgkelement &= \hdggelement  - \hdgfluxtouq_e^T \hdgaelement^{-1}\hdgfluxtouq_e\label{eq:hdg_operator} \eqdot
                 \intertext{The matrices~$\hdggelement$ and~$\hdgfluxtouq_e^T$ constitute the last three rows of the element operator in~\eqref{eq:ldg_element_operator_threed}:}
  \hdgfluxtouq_e^T &=
                     \begin{pmatrix}
                       \hdgfluxtou_{1,e}^T &  \hdgfluxtoq_{1,e}^T & \zeromat          & \zeromat          \\
                       \hdgfluxtou_{2,e}^T &  \zeromat          & \hdgfluxtoq_{2,e}^T & \zeromat          \\
                       \hdgfluxtou_{3,e}^T &  \zeromat          & \zeromat          & \hdgfluxtoq_{3,e}^T          \\
                     \end{pmatrix}\\
                 \hdggelement &=
                 \begin{pmatrix}
                   \hdgg_{1,e} & \zeromat   & \zeromat \\
                   \zeromat   & \hdgg_{2,e} & \zeromat \\
                   \zeromat   & \zeromat    & \hdgg_{3,e}
                 \end{pmatrix}\eqdot
  \intertext{%
  Here, $\hdggelement$ is block diagonal and integrates separately over each of the six faces and~$\hdgfluxtouq_e$ computes the effect of the fluxes~$\vecutilde_{i}$ onto~$\vecu$ and each~$\vecq_{i}$ in the element.
  The matrix~$\hdgaelement$ couples~$\vecu$ and~$\vecqvec$ in the element and is an extension of~\eqref{eq:element_operator_u_q_1d}
  }
    \hdgaelement &=
                   \begin{pmatrix}
                     \lambda \massmat_e + \hdgeelement & - \hdgdiffmat_{1,e}  & - \hdgdiffmat_{2,e}  & - \hdgdiffmat_{3,e}            \\
                     - \hdgdiffmat_{1,e}^T    & -\massmat_e  & \zeromat   & \zeromat   \\
                     - \hdgdiffmat_{2,e}^T    & \zeromat   & -\massmat_e & \zeromat  \\
                     - \hdgdiffmat_{3,e}^T    & \zeromat & \zeromat& -\massmat_e  \\
                   \end{pmatrix} \eqdot
\end{align}
As before, a~\schur\ complement allows to compute an operator coupling only~$\vecu$ with itself which takes the form
\begin{align}
  &\begin{aligned}
        \hdgzelementthreed^{-1} =
        \helmholtzparameter \geomcoeff{0,e}&\tp{\massmatstandard}{\massmatstandard}{\massmatstandard}
        + \geomcoeff{1,e} \tp{\massmatstandard}{\massmatstandard}{\stiffmat}\\
        + \geomcoeff{2,e} &\tp{\massmatstandard}{\stiffmat}{\massmatstandard}
        + \geomcoeff{3,e} \tp{\stiffmat}{\massmatstandard}{\massmatstandard}     \eqcomma
  \end{aligned}\label{eq:element_operator_u_q_3d}
  \intertext{leading to the inverse}
                                           &\hdgaelement^{-1} =
                                             \hdguqtou_e^T\hdgzelementthreed\hdguqtou_e
                                                 -
          \begin{pmatrix}
            \zeromat & \zeromat & \zeromat & \zeromat\\
            \zeromat & \massmat_e^{-1} & \zeromat & \zeromat \\
            \zeromat & \zeromat & \massmat_e^{-1} & \zeromat\\
            \zeromat & \zeromat &  \zeromat & \massmat_e^{-1}
          \end{pmatrix}  \label{eq:element_operator_u_3d} \eqcomma
                                             \intertext{where}
                  &\hdguqtou_e =  \begin{pmatrix}
                     \tp{\identitymat}{\identitymat}{\identitymat}\\
                     -\frac{2}{h_{1,e}}\tp{\identitymat}{\identitymat}{(\massmatstandard^{-1}\hdgdiffmatstandard^T)}\\
                     -\frac{2}{h_{2,e}}\tp{\identitymat}{(\massmatstandard^{-1}\hdgdiffmatstandard^T)}{\identitymat}\\
                     -\frac{2}{h_{3,e}}\tp{(\massmatstandard^{-1}\hdgdiffmatstandard^T)}{\identitymat}{\identitymat}\\
                   \end{pmatrix}^T \eqdot
\end{align}
The combination of~\eqref{eq:hdg_operator}, \eqref{eq:element_operator_u_q_3d}, and~\eqref{eq:element_operator_u_3d} allows for a tensor-product evaluation of the~HDG element operator, albeit one scaling with~$\order{\poly^4}$.

\section{A linearly scaling~\hdg\ operator in three dimensions}\label{sec:linearization}

\subsection{Sum factorization of the~\hdg\ operator}

Using the operators from~\prettyref{tab:hdg_matrices_threed}, all suboperators occurring in~\eqref{eq:hdg_operator} can be expressed as tensor products.
Similarly to~\eqref{eq:hdg_operator_oned}, two operations are present:
The matrix~$\hdggelement$ integrates the fluxes on each element boundary and is block-diagonal and, for a basis with a non-diagonal mass matrix, can be implemented in~$\order{\npoint^3}$ multiplications, whereas a diagonal mass matrix streamlines it to~$\order{\npoint^2}$.
Therefore, this part of the operator scales linearly already.
The hybridized term, however, requires a closer look.

Typically, the operator gets evaluated by first applying~$\hdgfluxtouq_e$, then~$\hdgaelement^{-1}$ and, lastly, mapping back~\cite{yakovlev_2015_hdg}.
The application of~$\hdgaelement^{-1}$ is often sped up using an extension of the generalized eigenvalue decomposition~\eqref{eq:generalized_eigenvalue_decomposition}: the fast diagonalization technique from~\cite{deville_2002_sem, lynch_1964_tensors}
\begin{subequations}
  \begin{align}
    \hdgzelementthreed &= \ptp{\extransmat}{\extransmat}{\extransmat} \diagonal_{\hdgzelementthreed}^{-1}\underbrace{\ptp{\extransmat^T}{\extransmat^T}{\extransmat^T}}_{\extransmat_e^T}\\
        \diagonal_{\hdgzelementthreed} &=
      \helmholtzparameter \geomcoeff{0,e}\tp{\identitymat}{\identitymat}{\identitymat}
      + \geomcoeff{1,e}\tp{\identitymat}{\identitymat}{\eigenvaluemat}
      + \geomcoeff{2,e}\tp{\identitymat}{\eigenvaluemat}{\identitymat}
      + \geomcoeff{3,e}\tp{\eigenvaluemat}{\identitymat}{\identitymat} \eqcomma
  \end{align}%
\end{subequations}
where~$\extransmat$ and~$\eigenvaluemat$ are the same as in the one-dimensional case, i.e.\,~\eqref{eq:generalized_eigenvalue_decomposition}.
However, the application of the three-dimensional tensor product~$\tp{\extransmat}{\extransmat}{\extransmat}$ still requires~$\order{\poly^4}$ operations and, hence, must be eliminated to lower the operation count to~${\order{\ndof} = \order{\poly^3\nelement}}$.
Here, a similar strategy as in~\cite{huismann_2014_condensation} proved sufficient:
Instead of using the operators~$\hdgfluxtouq_e$, $\hdguqtou_e$, and~$\tp{\extransmat}{\extransmat}{\extransmat}$ one after the other, they are fused together using~\eqref{eq:tp_property_fusing}.
This results in an operator mapping directly from the faces into the element eigenspace, where~$\diagonal_{\hdgzelementthreed}^{-1}$  can be applied, and then mapping back from the eigenspace directly to the faces.
As in the one-dimensional case, the approach requires an explicit form of~$\hdgbssub{e} = \extransmat_e^T\hdguqtou_e\hdgfluxtouq_e$, which computes to
\begin{align}
  \hdgbssub{e} &=
                 \begin{pmatrix}
                   \geomcoeff{1,e}\tp{\extransmat^T\massmatstandard}{\extransmat^T\massmatstandard}{\hdgbs} &
                   \geomcoeff{2,e}\tp{\extransmat^T\massmatstandard}{\hdgbs}{\extransmat^T\massmatstandard} &
                   \geomcoeff{3,e}\tp{\hdgbs}{\extransmat^T\massmatstandard}{\extransmat^T\massmatstandard}
                 \end{pmatrix} \eqdot
\end{align}
Using the above representation, the~HDG element operator simplifies as follows
\begin{align}
  \hdgkelement &= \hdggelement -  \hdgfluxtouq_e^T \hdgaelement^{-1} \hdgfluxtouq_e\nonumber\\
 \Rightarrow  \hdgkelement &= \hdggelement +  \hdgfluxtouq_e^T
                 \begin{pmatrix}
                   \zeromat & \zeromat & \zeromat & \zeromat \\
                   \zeromat & \massmat_e^{-1} & \zeromat & \zeromat \\
                   \zeromat & \zeromat & \massmat_e^{-1} & \zeromat \\
                   \zeromat & \zeromat & \zeromat & \massmat_e^{-1}
                 \end{pmatrix}
                                                    \hdgfluxtouq_e
                                                    -  \hdgfluxtouq_e^T\hdguqtou_e^T \hdgzelement^{-1}\hdguqtou_e   \hdgfluxtouq_e
 \nonumber \\
  \Rightarrow \hdgkelement  &= \hdggelement +
                 \begin{pmatrix}
                   \hdgfluxtoq_{1,e}^T\massmat_e^{-1}\hdgfluxtoq_{1,e} &\zeromat&\zeromat\\
                   \zeromat & \hdgfluxtoq_{2,e}^T\massmat_e^{-1}\hdgfluxtoq_{2,e} &\zeromat\\
                   \zeromat&\zeromat & \hdgfluxtoq_{3,e}^T\massmat_e^{-1}\hdgfluxtoq_{3,e}\\
                 \end{pmatrix}
 -  \hdgbssub{e}^T\diagonal_{\hdgzelementthreed}^{-1}\hdgbssub{e} \eqdot \label{eq:hdg_operator_factorized}
\end{align}
The first two terms implement the interaction between opposing faces of the element, using~$\order{\poly^3}$ multiplications for a non-diagonal mass matrix and~$\order{\poly^2}$ for a diagonal one.
The last term applies~$\hdgbssub{e}$, with~$\order{\poly^3}$ multiplications when using~$\tptwo{\extransmat^T\massmatstandard}{\extransmat^T\massmatstandard}$ on the faces and then expanding into the element eigenspace, then~$\diagonal_{\hdgzelementthreed}^{-1}$ in~$(\poly+1)^3$ multiplications and, lastly, maps back with~$\hdgbssub{e}^T$, also using~$\order{\poly^3}$ multiplications when reducing to the faces first.
Therefore, the operator can be applied in linear runtime, achieving linear scaling and, therefore, a major goal of the paper.

\begin{algorithm}[t]
  \caption{%
    Computation of the effect of the HDG operator on a per-element basis in the three-dimensional case, called~$\vecRtilde$, from the current flux~$\vecutilde$.
    For the tensor-product structure, the vectors~$\vecutilde_{i,  e}$ have the size~${2\times(\poly+1)\times(\poly+1)}$, ${(\poly+1)\times2\times(\poly+1)}$, and~${(\poly+1)\times(\poly+1)\times2}$ for the~$\coordinatesone$, $\coordinatestwo$, and~$\coordinatesthree$ direction, respectively.
  }%
  \label{alg:hdg_op_3d}
  \begin{algorithmic}
    \Function{HDG\_Op\_3D}{$\vecutilde$}
    \For{$ e = 1, \nelement$}
    \State{$\vecF_{\eigenspace} \gets \phantom{+} \ptp{\extransmat^T\massmatstandard}{\extransmat^T\massmatstandard}{\hdgbs}\geomcoeff{1,e}\vecutilde_{1,  e} $} \Comment[.3 \linewidth]{$\coordinatesone$ contribution to~$\vecF_{\eigenspace}$}
    \State{$\phantom{\vecF_{\eigenspace} \gets } + \ptp{\extransmat^T\massmatstandard}{\hdgbs}{\extransmat^T\massmatstandard}\geomcoeff{2,e}\vecutilde_{2,  e} $} \Comment[.3 \linewidth]{$\coordinatestwo$ contribution to~$\vecF_{\eigenspace}$}
    \State{$\phantom{\vecF_{\eigenspace} \gets } + \ptp{\hdgbs}{\extransmat^T\massmatstandard}{\extransmat^T\massmatstandard}\geomcoeff{3,e}\vecutilde_{3,  e} $} \Comment[.3 \linewidth]{$\coordinatesthree$ contribution to~$\vecF_{\eigenspace}$}
    \State{$\vecu_{\eigenspace} \gets \diagonal_{\hdgzelementthreed}^{-1} \vecF_{\eigenspace}$} \Comment[.3 \linewidth]{solution in eigenspace}
    \State{$\vecRtilde_{1,  e}    \gets \geomcoeff{1,e}\tp{\massmatstandard}{\massmatstandard}{(\hdgg  + \hdgfluxtoq^T \massmatstandard^{-1}\hdgfluxtoq)} \ \vecutilde_{1,  e} - \geomcoeff{1,e} \ptp{\massmat\extransmat}{\massmat\extransmat}{\hdgbs^T} \vecu_{\eigenspace}$}
    \State{$\vecRtilde_{2,  e}    \gets \geomcoeff{2,e}\tp{\massmatstandard}{(\hdgg  + \hdgfluxtoq^T \massmatstandard^{-1}\hdgfluxtoq)}{\massmatstandard} \ \vecutilde_{2,  e} - \geomcoeff{2,e} \ptp{\massmat\extransmat}{\hdgbs^T}{\massmat\extransmat} \vecu_{\eigenspace}$}
    \State{$\vecRtilde_{3,  e}    \gets \geomcoeff{3,e}\tp{(\hdgg  + \hdgfluxtoq^T \massmatstandard^{-1}\hdgfluxtoq)}{\massmatstandard}{\massmatstandard} \ \vecutilde_{3,  e} - \geomcoeff{3,e} \ptp{\hdgbs^T}{\massmat\extransmat}{\massmat\extransmat} \vecu_{\eigenspace}$}
    \EndFor
    \Return $\vecRtilde$
    \EndFunction
  \end{algorithmic}
\end{algorithm}
\prettyref{alg:hdg_op_3d} depicts an implementation of the operator which expands~\prettyref{alg:hdg_op_1d} to three-dimensional cuboidal tensor-product elements.
In every element, the six faces infer a residual for~$\vecu$ in the eigenspace~$\eigenspace$, called~$\vecF_{\eigenspace}$.
There, the inverse eigenvalues are applied to compute the inferred solution~$\vecu_{\eigenspace}$.
Lastly, the effect of this eigenspace solution and the interaction between the opposing faces is computed, resulting in the overall residual~$\vecRtilde_{e}$.

Computing the right-hand side in the element eigenspace in~\prettyref{alg:hdg_op_3d} requires~${24 \npoint^3}$ floating point operations when applying~$\extransmat^T\massmatstandard$, and further~${12 \npoint^3}$ operations for mapping into the eigenspace.
The application of~$\diagonal_{\hdgzelementthreed}^{-1}$ incurs another~$\npoint^3$ operations, and mapping back leads to another~$36 \npoint^3$.
Compared to these~$73\npoint^3$ operations, the two other terms remain insignificant:
Assuming a diagonal mass matrix, as present for~\legendre\ polynomials, \lagrange\ polynomials on~\gauss\ points, or when approximating the mass matrix for~\lagrange\ polynomials on~\gausslobatto\ nodes, the first two terms in~\eqref{eq:hdg_operator_factorized} only require~$\order{\poly^2}$ operations.

\subsection{Product factorization of the HDG operator}\label{sec:product_factorization}

\prettyref{alg:hdg_op_3d} allows for an evaluation of the~HDG residual in linear runtime.
However, the application cost is high when compared to a tensor-product operator for the primal form of~DG~\cite{stiller_2017a_multigrid, kronbichler_2017_hdg}.
Of the~$73\npoint^3$ operations, most apply only the two transformations, either~$\tptwo{\extransmat^T\massmatstandard}{\extransmat^T\massmatstandard}$ or~$\tptwo{\massmatstandard\extransmat}{\massmatstandard\extransmat}$.
These can, however, be removed by transforming the operator.
Consider the element operator mapping from the first direction to the first direction
\begin{align}
  &\begin{aligned}
    \hdgk_{11,e}
    &=       \geomcoeff{1,e}\tp{\massmatstandard}{\massmatstandard}{(\hdgg+\hdgfluxtoq^T \massmatstandard^{-1}\hdgfluxtoq})\\
    & -\geomcoeff{1,e}\ptp{\massmatstandard\extransmat}{\massmatstandard\extransmat}{\hdgbs^T}\diagonal_{\hdgzelementthreed}^{-1}\ptp{\extransmat^T\massmatstandard}{\extransmat^T\massmatstandard}{\hdgbs}\geomcoeff{1,e} \eqdot
  \end{aligned}
      \intertext{Applying~$\tp{\extransmat^T}{\extransmat^T}{\identitymat}$ from the left and~$\tp{\extransmat}{\extransmat}{\identitymat}$ from the right leads to}
    &    \hdgkhat_{11,e} = \ptp{\extransmat^T}{\extransmat^T}{\identitymat^T}\ \hdgk_{11,e}\ \ptp{\extransmat}{\extransmat}{\identitymat}\\
  &\begin{aligned}
    \hdgkhat_{11,e}
    &=     \ptp{\extransmat^T}{\extransmat^T}{\identitymat}\left[  \geomcoeff{1,e}\tp{\massmatstandard}{\massmatstandard}{(\hdgg+\hdgfluxtoq^T \massmatstandard^{-1}\hdgfluxtoq)} \right]\ptp{\extransmat}{\extransmat}{\identitymat}\\
    & -\geomcoeff{1,e}\left[\ptp{\extransmat^T}{\extransmat^T}{\identitymat}\ptp{\massmatstandard\extransmat}{\massmatstandard\extransmat}{\hdgbs^T}\diagonal_{\hdgzelementthreed}^{-1}\ptp{\extransmat^T\massmatstandard}{\extransmat^T\massmatstandard}{\hdgbs}\ptp{\extransmat}{\extransmat}{\identitymat} \right]\geomcoeff{1,e}
  \end{aligned}\\
    &\begin{aligned}
   \Rightarrow \hdgkhat_{11,e}
    &=     \geomcoeff{1,e}\left[\tp{\extransmat^T\massmatstandard\extransmat}{\extransmat^T\massmatstandard\extransmat}{({\hdgg+\hdgfluxtoq^T \massmatstandard^{-1}\hdgfluxtoq})}\right]\\
    & -\geomcoeff{1,e}\left[\ptp{\extransmat^T\massmatstandard\extransmat}{\extransmat^T\massmatstandard\extransmat}{\hdgbs^T}\diagonal_{\hdgzelementthreed}^{-1}\ptp{\extransmat^T\massmatstandard\extransmat}{\extransmat^T\massmatstandard\extransmat}{\hdgbs} \right]\geomcoeff{1,e}
  \end{aligned}
      \intertext{The identity~$\extransmat^T\massmatstandard\extransmat= \identitymat$ from~\eqref{eq:generalized_eigenvalue_decomposition_massmat} simplifies these terms to}
      &\hdgkhat_{11,e}
    =     \geomcoeff{1,e}\tp{\identitymat}{\identitymat}{({\hdgg+\hdgfluxtoq^T \massmatstandard^{-1}\hdgfluxtoq})}
     -\geomcoeff{1,e}\ptp{\identitymat}{\identitymat}{\hdgbs^T}\diagonal_{\hdgzelementthreed}^{-1}\ptp{\identitymat}{\identitymat}{\hdgbs} \geomcoeff{1,e} \eqdot
\end{align}
The terms for the faces in the~$x_2$ and~$x_3$ direction can be treated similarly by permutation of~$\tp{\extransmat^T}{\extransmat^T}{\identitymat}$, as can terms coupling different directions.
Therefore, every occurrence of~${\extransmat^T\massmatstandard}$ in the operator becomes~${\extransmat^T\massmatstandard\extransmat = \identitymat}$, and similarly~${\massmatstandard\extransmat}$ and~$\massmatstandard$ transform to identity as well.
While the transformation can simplify the implementation of~\prettyref{alg:hdg_op_3d}, applying the transformation to solution variable and residual instead lowers the operator costs across the whole solution process.

\begin{algorithm}[t]
  \caption{%
    Computation of the effect of the HDG operator on a per-element basis for the transformed basis in the three-dimensional case, called~$\vecRtildehat$, from the transformed flux~$\vecutildehat$.
  }%
  \label{alg:hdg_op_3d_transformed}
  \begin{algorithmic}
    \Function{HDG\_Op\_3D}{$\vecutildehat$}
    \For{$ e = 1, \nelement$}%
    \State{%
      $\vecF_{\eigenspace} \gets \ptp{\identitymat}{\identitymat}{\hdgbs}\geomcoeff{1,e}\vecutildehat_{1,  e}
      +   \ptp{\identitymat}{\hdgbs}{\identitymat}\geomcoeff{2,e}\vecutildehat_{2,  e}  + \ptp{\hdgbs}{\identitymat}{\identitymat}\geomcoeff{3,e}\vecutildehat_{3,  e}$}
    \State{$\vecu_{\eigenspace} \gets \diagonal_{\hdgzelementthreed}^{-1} \vecF_{\eigenspace}$}
    \State{$\vecRtildehat_{1,  e}    \gets \geomcoeff{1,e}\ \tp{\identitymat}{\identitymat}{(\hdgg   + \hdgfluxtoq^T \massmatstandard^{-1}\hdgfluxtoq)}\  \vecutildehat_{1,  e} - \geomcoeff{1,e} \ptp{\identitymat}{\identitymat}{\hdgbs^T} \vecu_{\eigenspace}$}
    \State{$\vecRtildehat_{2,  e}    \gets \geomcoeff{2,e}\ \tp{\identitymat}{(\hdgg + \hdgfluxtoq^T \massmatstandard^{-1}\hdgfluxtoq)}{\identitymat}\  \vecutildehat_{2,  e} - \geomcoeff{2,e} \ptp{\identitymat}{\hdgbs^T}{\identitymat} \vecu_{\eigenspace}$}
    \State{$\vecRtildehat_{3,  e}    \gets \geomcoeff{3,e}\ \tp{(\hdgg + \hdgfluxtoq^T \massmatstandard^{-1}\hdgfluxtoq)}{\identitymat}{\identitymat}\  \vecutildehat_{3,  e} - \geomcoeff{3,e} \ptp{\hdgbs^T}{\identitymat}{\identitymat} \vecu_{\eigenspace}$}
    \EndFor
    \Return $\vecRtilde$
    \EndFunction
  \end{algorithmic}
\end{algorithm}
\prettyref{alg:hdg_op_3d_transformed} depicts the resulting algorithm.
Not only does the operator now mainly consist of the mapping from the transformed faces into the element eigenspace and back, requiring~$25\npoint^3$ multiplications, the terms coupling opposing faces now use~$\order{\poly^2}$ in every basis, independent of the mass matrix.
While the algorithm requires transformation of the coefficients for the flux~${\vecutilde \to \vecutildehat}$ prior to the application, this can be done after computing the right-hand side for fluxes and, therefore, as a tensor-product operation on the faces, requiring~$\order{\poly^3}$ operations.

\subsection{Operator runtimes}\label{sec:operator_runtimes}

While the multiplication count constitutes a good measure for the asymptotic performance of an algorithm, it is typically insufficient to predict the efficiency attained in practice.
For instance, matrix-matrix multiplications require more operations than tensor-products, but can remain faster over a wide range of polynomial degrees~\cite{cantwell_2011_operator, huismann_2018_optimization}.
Therefore, runtime tests were conducted to validate both, that the linear scaling is achieved as well as that matrix-matrix implementations are outperformed.

Two variants of the~HDG element operator were considered.
The first implemented~\prettyref{alg:hdg_op_3d} for~\gllnames~(GLL) polynomials.
As it implements a~HDG operator using tensor products it was called~HDG-TP\@.
The second one implemented~\prettyref{alg:hdg_op_3d_transformed}, which utilizes a transformation to streamline the operation count, and was therefore named~HDG-TPT.
Variants for the statically condensed continuous Galerkin case~(CCG) implementing the static condensed operator served as reference.
The first thereof implemented the condensed part using one matrix-matrix multiplication and was, therefore, called~CCG-MM\@.
The two other variants used a tensor-product decomposition, the first one on~GLL polynomials, called~CCG-TP, and the second one, CCG-TPT, using a transformed coordinate system and thereby lowering the number of operations.
Derivations of these can be found in~\cite{huismann_2016_condensation, huismann_2017_condensation}.
The matrix-matrix-based operator required~$\order{\poly^4}$ operations, whereas~$\order{\poly^3}$ operations sufficed for all tensor-product-based operators.
\prettyref{tab:operator_variants} lists the operators, leading factors of the operation counts, and associated preprocessing costs.

\begin{table}[t]
  \centering
  \caption{Application and preprocessing costs for the operators in conjunction with the incurred number of loads and stores per operator evaluation.
    Here, CCG denotes the condensed continuous~\galerkin\ method, whereas~HDG denotes the hybridized discontinuous~\galerkin\ method.
    For the condensed case, a factor of three was imposed on the loads and stores of to account for the evaluation of the primary part.
  }%
  \label{tab:operator_variants}
  \begin{tabular}{llllllr}\toprule
    Operator && Preprocessing cost && \# FLOP && \# loads / stores \\\midrule
    CCG-MM  && $\order{\poly^5}$ && $72 (\poly-1)^4 + \order{\poly^2}$ && $3 \cdot 18 {(\poly-1)}^2 \phantom{\ + {(\poly-1)}^3}$ \\
    CCG-TP  && $\order{\poly^3}$ && $97 (\poly-1)^3 + \order{\poly^2}$ && $3 \cdot 18 {(\poly-1)}^2 + {(\poly-1)}^3$ \\
    CCG-TPT && $\order{\poly^3}$ && $25 (\poly-1)^3 + \order{\poly^2}$ && $3 \cdot 18 {(\poly-1)}^2 + {(\poly-1)}^3$ \\
    HDG-TP  && $\order{\poly^3}$ && $73 (\poly+1)^3 + \order{\poly^2}$ && $        18 {(\poly+1)}^2 + {(\poly+1)}^3$ \\
    HDG-TPT && $\order{\poly^3}$ && $25 (\poly+1)^3 + \order{\poly^2}$ && $        18 {(\poly+1)}^2 + {(\poly+1)}^3$           \\\bottomrule
  \end{tabular}
\end{table}

The operators were implemented in Fortran 2008 using double precision.
The Intel Fortran compiler~\intelcompilerversion\ compiled the variants, with the associated~MKL serving as BLAS implementation.
A single core of an Intel~Xeon~E5-2680 v.3 constituted the measuring platform.
As it ran with a clock speed of~\SI{2.5}{\giga\hertz}, it allowed for a maximum floating point rate of~\SI{40}{\giga\flops}~\cite{hackenberg_2015_energy}.

For every polynomial degree in~$\{2\dots 32\}$, the operators were used~$101$ times on~${\nelement = 8^3}$ spectral elements, with the runtime of the last~$100$ times being measured with~\texttt{MPI\_Wtime}.
This approach precludes measurement of instantiation effects, for instance of libraries such as~BLAS\@.

\begin{figure}[t]
  \hspace*{\fill}
  \includegraphics{hdg_op_plot_condensed_operator_comparison_prep_legend.pgf}
  \hspace*{\fill}\\
  \hspace*{\fill}
  \begin{subfigure}{0.49\linewidth}
    \includegraphics{hdg_op_plot_condensed_operator_comparison_prep_times.pgf}
    \caption{}\label{subfig:op_runtimes_prep_times}
  \end{subfigure}
  \hfill
  \begin{subfigure}{0.49\linewidth}
    \includegraphics{hdg_op_plot_condensed_operator_comparison_comp_times.pgf}
    \caption{}\label{subfig:op_runtimes_comp_times}
  \end{subfigure}
  \hspace*{\fill}\\
  \hspace*{\fill}
  \begin{subfigure}{0.49\linewidth}
  \includegraphics{hdg_op_plot_condensed_operator_comparison_dof_per_s.pgf}
    \caption{}\label{subfig:op_runtimes_dof_per_s}
  \end{subfigure}
  \hfill
  \begin{subfigure}{0.49\linewidth}
  \includegraphics{hdg_op_plot_condensed_operator_comparison_gflops.pgf}
    \caption{}\label{subfig:op_runtimes_gflops}
  \end{subfigure}
  \hspace*{\fill}
  \caption{Operator runtimes when varying the polynomial degree~$\poly$ in a homogeneous mesh consisting of~${\nelement = 8^3 = 512}$ spectral elements.
    \subref{subfig:op_runtimes_prep_times}:~Preparation times for the matrices used by the operators.
    \subref{subfig:op_runtimes_comp_times}:~Operator runtimes.
    \subref{subfig:op_runtimes_dof_per_s}:~Rate of updates per second, where the equivalent number of degrees of freedom is calculated to be~${(\poly+1)}^3\nelement$.
    \subref{subfig:op_runtimes_gflops}:~Rate of floating point operations per second measured in~\SI{}{\giga\flops}.}%
  \label{fig:operator_runtimes}
\end{figure}
\prettyref{fig:operator_runtimes} depicts the measured instantiation times and runtimes in combination with the achieved floating point rate and equivalent degrees of freedom per second.
For the preparation time, the expected orders are achieved:
Computing the assembled element matrix for the static condensed method~CCG-MM scales with~$\order{\poly^5}$.
The preprocessing time of the other operators scales with~$\order{\poly^3}$, which can be attributed to them only requiring solution of an eigenvalue problem and then computing the inverse eigenvalues, both scaling with~$\order{\poly^3\nelement}$.

The runtimes of the operators, shown in~\prettyref{subfig:op_runtimes_comp_times}, fall into three categories.
The first one consists of the matrix-matrix-based version~CCG-MM\@.
While attaining the highest rate of floating point operations, near~\SI{35}{\giga\flops} beyond~$\poly=8$, the runtime increases significantly compared to the other methods, stemming from an increasing number of required operations per degree of freedom.
This is reflected in the number of updated degrees of freedom which decreases substantially with~$\poly$.
The non-transformed operators~HDG-TP and~CCG-TP constitute the second category, generating a mostly constant rate of updates while exhibiting a smaller runtime than~CCG-MM starting from~${\poly=6}$.
For large polynomial degrees the rate of updates oscillates depending on the polynomial degree.
The maximum is achieved when the number of points per direction~$\poly+1$ is a multiple of four, which lends itself to optimization for the given architecture, whereas slightly lower update rates result at other polynomial degrees.
For~CCG-TP the operator consists of a primary and a condensed part.
The latter resembles the~HDG operator and can be optimized well, but the former couples the vertices, edges, and faces of the element and consists of multiple smaller suboperators.
While oscillations occurs, they are not as large as for~HDG-TP\@.

Transforming the equation system allows leveraging the operators~HDG-TPT and~CCG-TPT\@.
These are faster than~CCG-MM for all polynomial degrees and gain a factor of four over the non-transformed variants.
The~HDG variant is slightly faster at lower polynomial degrees, due to the primary part of the condensed system being evaluated for~CCG-TPT\@.
For~${\poly > 4}$, however, both variants attain similar update rates, only at different polynomial degrees.
With the tensor-product variants exhibiting a mostly constant rate of updates, a major goal has been achieved: operators scaling linearly with the number of degrees of freedom.
Moreover, these attain the same degree of efficiency as the ones for the continuous discretization.

\begin{figure}
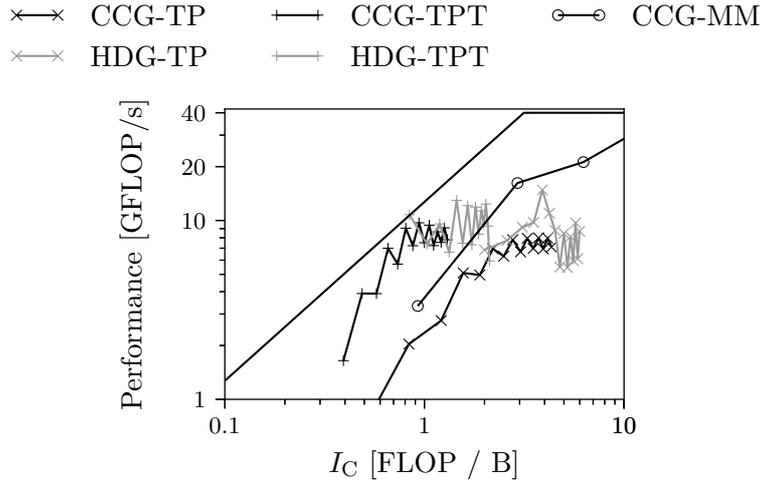

  \hspace*{\fill}
  \includegraphics{hdg_op_plot_condensed_operator_roofline_legend.pgf}
  \hspace*{\fill}\\
  \hspace*{\fill}
  \includegraphics{hdg_op_plot_condensed_operator_roofline.pgf}
  \hspace*{\fill}
  \caption{%
    Roofline analysis of the operator variants, depicting the performance attained in practice over the computational intensity~$I_{\mathrm{C}}$ resulting from~\prettyref{tab:operator_variants}.
    In the analysis, a limiting memory bandwidth of~\SI{12.7}{\giga\byte\per\second} and a maximum performance of~\SI{40}{\giga\flop\per\second} are used, as measured in~\cite{hackenberg_2015_energy}.
    Furthermore, for clarity of presentation, only odd polynomial degrees are shown, as these result in even operator sizes.
  }%
  \label{fig:operator_roofline}
\end{figure}%
Lastly, \prettyref{subfig:op_runtimes_gflops} depicts the rate of floating-point operations.
All variants start with a relatively low rate.
However, the rate increases quickly for the matrix-matrix variant, which attains more than~\SI{15}{\giga\flops} at~$\poly=5$ and saturates at~\SI{35}{\giga\flops}, which is beyond the range shown here.
These operations, however, are largely wasted since the number of updates decreases, as shown in~\prettyref{subfig:op_runtimes_dof_per_s}.
Due to evaluating the primary part, the variants~CCG-TP and~CCG-TPT for the static condensed operator harness only one eighth of the compute power at~$\poly=12$, whereas the variants for~HDG quickly use one quarter.
Beyond~$\poly=16$, the~HDG variants exhibit a highly oscillatory behaviour compared to the variants for static condensation.
This is due to the very short, monolithic operators.
These are far easier to optimize for the compiler, leading to spikes where the operator width is a multiple of the register width.
All tensor-product operators attain only near a quarter of the maximum rate.
This stems from the implementation, where the matrix of eigenvalues~$\diagonal_{\hdgzelementthreed}^{-1}$ is precomputed and stored explicitly.
In a roofline analysis, as shown in~\prettyref{fig:operator_roofline}, the usage of~$\diagonal_{\hdgzelementthreed}^{-1}$ leads to a constant as asymptotic limit for the computational intensity~\cite{williams_2009_performance}.
However, while loading~$\order{\poly^3}$ floating point numbers, the operation is still far from being memory-bound when computing at high polynomial degrees.
For machines with a lower critical computational intensity, computing the inverse eigenvalues in the operator itself can be a remedy, but was not found to be beneficial with the CPUs utilized in this study, even when capitalizing on imprecise floating point division.
However, the growing memory gap might render it beneficial in the future.

\section{Construction of an elliptic solver with linear runtime}\label{sec:runtimes}

\subsection{Linearly scaling face-local preconditioners}\label{sec:preconditioners}

With a linearly scaling residual evaluation, the door to linearly scaling solvers is wide open.
The only remaining obstacle is attaining a well-performing linearly scaling preconditioner.
This requirement severely limits the possibilities.
Assuming a constant iteration count, only tensor-product preconditioners working on the faces remain a possibility.
In general, a multigrid approach leads to very efficient solvers, as done for the continuous discretization in~\cite{haupt_2013_multigrid, huismann_2018_multigrid, huismann_2018_diss}.
There, \schwarz-type preconditioners using tensor-products were employed to generate an efficient smoother with linear scaling.
The derivation of a linearly-scaling~\schwarz-type smoother required for the multigrid algorithm is, however, a topic of its own and well beyond the scope of this paper.

To investigate the performance achievable with the linearly scaling operator, block-\jacobi\ methods are used.
They proved highly efficient for solution of diffusion equations in semi-implicit time-stepping schemes~\cite{huismann_2018_diss} and provide an astonishing resilience against high aspect ratios~\cite{huismann_2017_condensation}.
Moreover, they include the same kind of operations found in~\schwarz-type smoothers for multigrid and serve to evaluate their qualitative behaviour for high polynomial degrees.
Due to these positive properties the same approach is employed here as well.

Two preconditioners are derived:
The first one, applies the exact inverse on every face of the mesh, leading to a block-\jacobi\ preconditioner, whereas the second one only utilizes the inverse of the main diagonal.
To derive the preconditioners, a single element~$\elementdomain$ is considered.
The operator connecting two opposing faces, as illustrated here for the first direction, reads
\begin{align}
  &\begin{aligned}
    \hdgk_{11,e}
    &=       \geomcoeff{1,e}\tp{\massmatstandard}{\massmatstandard}{(\hdgg+\hdgfluxtoq^T \massmatstandard^{-1}\hdgfluxtoq})\\
    & -\geomcoeff{1,e}\ptp{\massmatstandard\extransmat}{\massmatstandard\extransmat}{\hdgbs^T}\diagonal_{\hdgzelementthreed}^{-1}\ptp{\extransmat^T\massmatstandard}{\extransmat^T\massmatstandard}{\hdgbs}\geomcoeff{1,e} \eqdot
  \end{aligned}\label{eq:hdg_op_east_east}
      \intertext{When removing the occurences of~$\massmatstandard\extransmat$ by multiplying with~$\extransmat^T$ from the left and~$\extransmat$ from the right in both face directions, the operator becomes}
    &\begin{aligned}
      \hdgkhat_{11,e} &=
      \ptp{\extransmat^T}{\extransmat^T}{\identitymat} \hdgk_{11, e} \ptp{\extransmat}{\extransmat}{\identitymat}\\
      &=       \geomcoeff{1,e}(\tp{\identitymat}{\identitymat}{\hdgg + \hdgbs^T \massmatstandard^{-1}\hdgbs)}
      -\geomcoeff{1,e}\ptp{\identitymat}{\identitymat}{\hdgbs^T}\diagonal_{\hdgzelementthreed}^{-1}\ptp{\identitymat}{\identitymat}{\hdgbs} \geomcoeff{1,e}\eqcomma
    \end{aligned}\label{eq:hdg_op_east_east_transformed}
\end{align}
which is diagonal in the two directions on the faces.
As each face occurs in two elements, two elements contribute to the preconditioner.
Let~$\mymat{Y}_e$ denote the effect of~$\hdgkelement$ from one face to itself for all faces.
Then, ${\mymat{\tilde{Y}} = \hdgscattermatrix\hdggathermatrix \mymat{Y}}$ is the effect of the global flux of the face onto itself.
However, the application requires the variable to reside in the face eigenspace, i.e.\,after applying the transformation matrix~$\extransmat^T\massmat$ in the two directions of the face.
Therefore, application of the preconditioner consists of first mapping into the face eigenspaces, applying the inverse~${\mymat{\tilde{Y}}^{-1}}$, and then mapping back.
\prettyref{alg:hdg_pc_block} shows an implementation thereof.
The preconditioner requires four one-dimensional matrix products on each face, resulting in~${4 \cdot 2 \cdot 6 \npoint^3 = 48\npoint^3\nelement}$ floating point operations per application of the preconditioner.
Moreover, as~$\hdgkhat_{ii,e}$ from~\eqref{eq:hdg_op_east_east_transformed} can be computed in~$\order{\poly^3}$ operations, determining~$\mymat{\tilde{Y}}^{-1}$ requires~$\order{\poly^3\nelement}$ operations.
\begin{algorithm}[t]
  \caption{%
    Application of the block preconditioner to the residual~$\vecRtilde$ in a non-transformed system.
    Here, all variables are stored and accessed in a face-wise fashion and~$\nfaces$ denote the total number of faces in the mesh.
  }%
  \label{alg:hdg_pc_block}
  \begin{algorithmic}
    \Function{Block\_Preconditioner}{$\vecRtilde$}
    \For{$ f = 1, \nfaces$} \Comment{Loop over all faces}
    \State{$\vecRtildehat_{f} \gets \ptptwo{\extransmat^T}{\extransmat^T}\vecRtilde_{f}$} \Comment{transformation}
    \State{$\vecztildehat_f \gets \mymat{\tilde{Y}}_{f}^{-1} \vecRtildehat_f$} \Comment{inversion on the faces}
    \State{$\vecztilde_{f} \gets \ptptwo{\extransmat}{\extransmat}\vecztildehat_{f}$} \Comment{transformation}
    \EndFor{}
    \Return{$\vecztilde$}
    \EndFunction{}
  \end{algorithmic}
\end{algorithm}

With~\prettyref{alg:hdg_pc_block}, an explicit representation for the inverse on the faces is present.
However, the application costs are nearly as high as the operator costs.
A diagonal preconditioner can solve the latter issue, at the expense of the condition of the system.
To compute the diagonal preconditioner for one point~$i,j$ on a face, the restrictor~${\mymat{e}_{j}\otimes\mymat{e}_{i}}$ can be utilized.
The diagonal results via~$({\mymat{e}_{j}\extransmat}\otimes{\mymat{e}_{i}\extransmat})\mymat{\tilde{Y}}^{-1}({\extransmat^T \mymat{e}_{j}^T}\otimes{\extransmat^T\mymat{e}_{i}^T})$, which can be precomputed in~$\order{\poly^3}$ and applied in~$\order{\poly^2}$.

\subsection{Summary of solution method}

So far only pieces of the solution process were shown.
This subsection states the solution method as a whole.
The solution process with~HDG works as follows:
First, hybridize the system by computing the fluxes~$\vecutilde$ from the initial guess for the solution variable~$\vecu$, first, computing~$\vecq$ via~\eqref{eq:helmholtz_flux_q} and then using~\eqref{eq:fluxu}.
For a non-primal formulation, which includes the auxiliary variable~$\vecq$ as solution variable, $\vecutilde$ can be computed directly via~\eqref{eq:fluxu}.
Thereafter, the right-hand side~$\hdgrhs$ is computed from the right-hand side for~$\vecu$ and the boundary conditions.
Solution of the~HDG system~\eqref{eq:hdg_system} constitutes the third step.
Lastly, the solution is recomputed from the flux.
\prettyref{alg:hdg_cg} summarizes the solution algorithm.

\begin{algorithm}[t]
  \caption{%
    Solution process with~HDG for a primal method.
  }%
  \label{alg:hdg_cg}
  \begin{algorithmic}
    \Function{HDG\_Solver}{$\vecu$, $\vecf$}
    \For{$ e = 1, \nelement$}%
    \State{$\vecq_{1,e} \gets \geomcoeff{1,e} \ptp{\identitymat}{\identitymat}{\diffmat} \vecu_{e} $}\Comment{Compute~$\vecq_{1}$ via differentiation}
    \State{$\vecq_{2,e} \gets \geomcoeff{2,e} \ptp{\identitymat}{\diffmat}{\identitymat} \vecu_{e} $}\Comment{Compute~$\vecq_{2}$ via differentiation}
    \State{$\vecq_{3,e} \gets \geomcoeff{3,e} \ptp{\diffmat}{\identitymat}{\identitymat} \vecu_{e} $}\Comment{Compute~$\vecq_{3}$ via differentiation}
    \State{$\hdgrhs_{e} \gets \mymat{g}_{\mathrm{N},e} + \hdgfluxtouq_e^T\hdgaelement^{-1}\begin{pmatrix} \massmat_{e}\vecf_{e} & \mymat{0} \end{pmatrix}^T$}
    \EndFor{}
    \State{$\vecutilde \gets \text{Flux}(\vecu, \vecq)$} \Comment{initial guess}
    \State{}
    \State{$\text{pCG\_Solve}{(\hdgk \vecutilde =  \hdgrhs)}$} \Comment{solution with pCG method}
    \State{}
    \For{$ e = 1, \nelement$}%
    \State{$\vecF_e \gets \begin{pmatrix} \massmat_{e}\vecf_{e} & \mymat{0} \end{pmatrix}^T + \hdgfluxtouq_e \vecutilde_e $} \Comment{Right-hand side for~$\vecu$ and~$\vecq$}
    \State{$\begin{pmatrix}\vecu_e& \vecq_e\end{pmatrix}^T \gets \hdga_e^{-1} \vecF_e$} \Comment{Recompute solution variables}
    \EndFor{}
    \Return{$\vecu$}
    \EndFunction{}
  \end{algorithmic}
\end{algorithm}%
In the solution process, the hybridization incurs a three-dimensional tensor product to transform~$\vecu$ into the element operator eigenspace.
This step scales with~$\order{\poly^4\nelement}$, as does the transformation of the right-hand side.
But at least these occur only once.
The solution process, when implemented with a~CG method, consists of an operator evalation followed by scalar products and a preconditioner application.
When using the operator evaluation from~\prettyref{sec:linearization} and preconditioners from~\prettyref{sec:preconditioners}, these are lowered in complexity from~$\order{\poly^4\nelement}$ to~${\order{\poly^3\nelement} = \order{\ndof}}$, i.e.\,the resulting solvers achieve a linearly scaling iteration time.
The last step again recquires a transformation into the element eigenspaces and, therefore, $\order{\poly^4\nelement}$ operations.
However, the pre- and postprocessing occur, by definition, once per solution process, whereas the iteration process is reiterated tens if not hundreds of times.
For the continuous method, the iteration, not pre- and postprocessing are the most costly component for polynomial degrees up to~${\poly > 48}$~\cite{huismann_2018_multigrid}, and these lie far outside the range of currently employed ones.
Therefore, the solver can be described as scaling linearly for all relevant polynomial degrees.

\begin{algorithm}[t]
  \caption{%
    Solution process with~HDG for the transformed system.
  }%
  \label{alg:hdg_transformed_cg}
  \begin{algorithmic}
    \Function{Transformed\_Solver}{$\vecutilde$, $\hdgrhs$}
    \For{$ f = 1, \nfaces$}
    \State{$\hdgrhshat_{f} \gets \ptptwo{\extransmat^T}{\extransmat^T}\hdgrhs_{f}$} \Comment{transformed right-hand side}
    \State{$\vecutildehat_{f} \gets \ptptwo{\extransmat^T\massmat}{\extransmat^T\massmat}\vecutilde_{f}$} \Comment{transformed initial guess}
    \EndFor{}
    \State{$\text{pCG\_Solve}{(\hdgkhat \vecutildehat =  \hdgrhshat)}$} \Comment{solution in transformed system}
    \For{$ f = 1, \nfaces$}
    \State{$\vecutilde_{f} \gets \ptptwo{\extransmat}{\extransmat}\vecutildehat_{f}$} \Comment{transformed initial guess}
    \EndFor{}
    \Return{$\vecutilde$}
    \EndFunction{}
  \end{algorithmic}
\end{algorithm}%
The described algorithm works as is for the untransformed system.
To leverage the product factorization from~\prettyref{sec:product_factorization}, a coordinate transformation is required, which gets prepended and appended to the solver call in~\prettyref{alg:hdg_cg}.
As both of these operations can work on the already hybridized data and can be written in tensor-product form, they scale linearly with the overall number of degrees of freedom and do not impede the solution process.
\prettyref{alg:hdg_transformed_cg} summarizes the altered solver.

For testing purposes, four~HDG solvers are considered, all leveraging the developed linear operators:
The first one is a direct~CG implementation of~HDG without preconditioning and called~\solverhdgunprec.
The second one, \solverhdgdiag, utilizes diagonal preconditioning.
The third one, \solverhdgblock, employs the block-preconditioner from~\prettyref{alg:hdg_pc_block}.
The last considered~HDG solver capitalizes on the transformation in~\prettyref{alg:hdg_transformed_cg}:
The application of operator and preconditioner both require application of~$\extransmat^T\massmat$ and~$\extransmat^T$, respectively, in the face directions before allowing for evaluation on the face.
By applying these as coordinate transformation, \prettyref{alg:hdg_op_3d_transformed} can be used directly for evaluation.
Moreover, the block preconditioner~\prettyref{alg:hdg_pc_block} simplifies to a diagonal one.
Hence, the transformation streamlines the operation count of both operator and block preconditioner, and the resulting solver is called~\solverhdgtrans.
For all of these solvers, a counterpart using the continuous discretization with static condensation is investigated as well.

\begin{table}[t]
  \centering
  \caption{Number of operations for precomputation of matrices~($N_{\mathrm{mat}}$), pre- and post-processing of variables~($N_{\mathrm{prepost}}$), application of operator~($N_{\mathrm{op}}$) and application of the preconditioner~($N_{\mathrm{cond}}$) for the proposed solvers.}%
  \label{tab:solver_multiplications}
  \begin{tabular}{llccccccrr}\toprule
    Solver && $N_{\mathrm{mat}}$ && $N_{\mathrm{prepost}}$ && $N_{\mathrm{op}}$ && $N_{\mathrm{cond}}$ \\\midrule
    \solverhdgunprec && $\order{\poly^3\nelement}$ && $\order{\poly^4\nelement}$ && $73 {(\poly+1)}^3\nelement$ && --  \\
    \solverhdgdiag   && $\order{\poly^3\nelement}$ && $\order{\poly^4\nelement}$ && $73 {(\poly+1)}^3\nelement$ && $6 {(\poly+1)}^2\nelement$  \\
    \solverhdgblock  && $\order{\poly^3\nelement}$ && $\order{\poly^4\nelement}$ && $73 {(\poly+1)}^3\nelement$ && $49 {(\poly+1)}^3\nelement$  \\
    \solverhdgtrans  && $\order{\poly^3\nelement}$ && $\order{\poly^4\nelement}$ && $27 {(\poly+1)}^3\nelement$ && $6 {(\poly+1)}^2\nelement$  \\\midrule
    \solverccgunprec && $\order{\poly^3\nelement}$ && $\order{\poly^4\nelement}$ && $97 {(\poly-1)}^3\nelement$ && --  \\
    \solverccgdiag   && $\order{\poly^3\nelement}$ && $\order{\poly^4\nelement}$ && $97 {(\poly-1)}^3\nelement$ && $6 {(\poly-1)}^2\nelement$  \\
    \solverccgblock  && $\order{\poly^3\nelement}$ && $\order{\poly^4\nelement}$ && $97 {(\poly-1)}^3\nelement$ && $49 {(\poly-1)}^3\nelement$  \\
    \solverccgtrans  && $\order{\poly^3\nelement}$ && $\order{\poly^4\nelement}$ && $25 {(\poly-1)}^3\nelement$ && $6 {(\poly-1)}^2\nelement$
    \\ \bottomrule
  \end{tabular}
\end{table}
\prettyref{tab:solver_multiplications} lists the pre- and post-processing times for the solvers.
For both, continuous and discontinuous discretization, the right-hand side results by applying an element-wise inverse via fast diagonalization.
Similarly, the generation of the solution consists by mapping into the element eigenspace and applying a tensor-product to the solution there.
These two operations scale with~$\order{\poly^4\nelement}$ and are only evaluated once during the solution process, such that they not impede the linear scaling.
They are not expected to dominate the runtime behaviour until at least~${\poly > 48}$~\cite{huismann_2018_multigrid}, allowing to disregard the issue here.

\subsection{Definition of a test case for assessment}\label{sec:solver_runtimes}

In the following the test case from~\cite{huismann_2017_condensation} serves to evaluate the solvers derived in the previous section.
In a domain~${\domain = {(0, 2\pi)}^3}$ the manufactured solution
\begin{align}
  &\begin{aligned}
    &u_{\mathrm{ex}}\of{x} = \cos\of{k (x_1 - 3 x_2 + 2 x_3)}   \sin\of{k (1 + x_1)}   \\
    &\cdot \sin\of{k (1 - x_2)}   \sin\of{k (2 x_1 + x_2)}   \sin\of{k (3 x_1 - 2 x_2 + 2 x_3)} \eqcomma
    \label{eq:solution}
  \end{aligned}
      \intertext{%
      is considered, with the stiffness parameter~$k$ set to 5, and the right-hand side of~\eqref{eq:helmholtz_equation} analytically evaluated from}
    &f \of{x} = \lambda u_{\mathrm{ex}}\of{x} - \Delta u_{\mathrm{ex}}\of{x}  \eqdot
\end{align}
Inhomogeneous~\dirichlet\ boundary conditions are used to impose the exact values, as defined by~\eqref{eq:solution}, on the boundary.
The initial guess for~$\vecu$ consists of pseudo-random numbers in the interior of the domain, with the~\dirichlet\ conditions imposed on the boundary.
This initial guess is used directly for the continuous discretization, whereas for~HDG the initial fluxes resulting from the strong form of~\eqref{eq:fluxu} are imposed.
While~${\helmholtzparameter > 0}$ was utilized in preliminary validations, the~\poisson\ case~${\lambda =  0}$ leads to a worse condition.
This generates a system that is harder to solver and is, therefore, utilized here to demonstrate the performance for the worst case.

As in~\prettyref{sec:operator_runtimes}, the solvers were implemented in Fortran using double precision, compiled with the~Intel Fortran compiler, and a single core of an Intel~Xeon~E5-2680 v3 constituted the measuring platform.
On this machine, the solvers were run~$11$ times.
The runtime of the last ten runs was measured with~\texttt{MPI\_Wtime}, whereas the first run precluded measurement of instantiation effects.
In every run, the iterations were stopped after the initial residual was reduced by a factor of~$10^{-10}$.

\subsection{Robustness against increase in polynomial degree}

The first test uses a constant number of elements~${\nelement = 8^3 = 512}$ while varying the polynomial degree in~${\poly \in \{2 \dots 32\} }$, with the penalty parameter set to~${\hdgtauelement = 25}$.

\begin{figure}[t]
  \hspace*{\fill}
  \includegraphics{hdg_op_plot_condensed_solver_po_legend.pgf}
  \hspace*{\fill}\\
  \begin{subfigure}{0.49\linewidth}
    \includegraphics{hdg_op_plot_condensed_solver_po_error_max.pgf}
    \caption{}\label{subfig:solver_runtimes_homogeneous_error_max}
  \end{subfigure}
  \begin{subfigure}{0.49\linewidth}
    \includegraphics{hdg_op_plot_condensed_solver_po_iters.pgf}
    \caption{}\label{subfig:solver_runtimes_homogeneous_iters}
  \end{subfigure}\\
  \hspace*{\fill}
  \begin{subfigure}{0.49\linewidth}
    \includegraphics{hdg_op_plot_condensed_solver_po_times.pgf}
    \caption{}\label{subfig:solver_runtimes_homogeneous_times}
  \end{subfigure}
  \begin{subfigure}{0.49\linewidth}
  \includegraphics{hdg_op_plot_condensed_solver_po_times_per_dof.pgf}
    \caption{}\label{subfig:solver_runtimes_homogeneous_times_per_dof}
  \end{subfigure}
  \caption{Results when varying the polynomial degree~$\poly$ in a homogeneous mesh consisting of~${\nelement = 8^3}$ spectral elements.
    For the~HDG solvers, the penalty parameter was fixed at~${\hdgtauelementi = 25}$.
    \subref{subfig:solver_runtimes_homogeneous_error_max}:~Resulting discretization error.
    \subref{subfig:solver_runtimes_homogeneous_iters}:~Number of iterations required to lower the residual by ten orders.
    \subref{subfig:solver_runtimes_homogeneous_times}:~Solver runtimes.
    \subref{subfig:solver_runtimes_homogeneous_times_per_dof}:~Runtimes per unknown.}%
  \label{fig:solver_runtimes_homogeneous}
\end{figure}
\prettyref{fig:solver_runtimes_homogeneous} depicts the solution error, the resulting number of iterations, the runtime, and the runtime per iteration and unknown.
The error of all computed solutions behaves as it should.
With the given number of elements and smaller polynomial degrees the discrete solution is not yet in the convergence range for this value of the stiffness parameter.
Beyond~${\poly \approx 12}$ spectral convergence is noticed, reflected by an increasing slope with increasing~$\poly$, validating the discretization and solution procedures.
Beyond a value of~${\poly\approx 26}$ to~${\poly\approx 30}$ the error saturates for all methods.
The latter behaviour is common with high-order methods and results from a limited machine accuracy in combination with an accumulation of round-off errors.
The level at which this occurs depends on the method.
Here, the discontinuous discretization exhibits a somewhat higher error level compared to the continuous one.
Furthermore, a difference between the unpreconditioned, and preconditioned variants occurs:
The former have a higher remaining error, which can stem from the combination of preconditioning in conjunction with a CG algorithm that does not include a restart, possibly leading to an accumulation of round-off errors and searches in the wrong subspace~\cite{shewchuk_1994_cg}.

When comparing the number of iterations in~\prettyref{subfig:solver_runtimes_homogeneous_iters}, HDG always requires more than the continuous formulation.
For a low polynomial degree, a factor of two is present, whereas continuous and discontinuous variants converge against each other for high polynomial degrees.
This indicates that the~HDG has a worse condition number, stemming from a combination of more degrees of freedom and the ill-conditioned penalty terms, which dominate until high polynomial degrees.

As to be expected, a diagonal preconditioner for HDG results in a lower required number of iterations compared to the unpreconditioned version.
The gain remains, however, small.
In contrast, the block-preconditioner leads to a significant improvement, with~$\approx 100$ iterations being required over a wide range of polynomial degrees, and a factor of two less iterations at~${\poly=32}$.
While the number of iterations remains higher than in the continuous case for~$\poly<16$, the difference vanishes when further increasing the order.

\prettyref{subfig:solver_runtimes_homogeneous_times} depicts the runtimes.
The unpreconditioned and diagonally-preconditioned solvers incur the highest runtime per degree of freedom, with the benefit of diagonal preconditioning mostly being offset by the higher costs.
The block-preconditioned versions in the non-transformed systems, in contrast, exhibit a constant runtime per degree of freedom for~$\poly\geq8$, which lies near~\SI{2.5}{\micro\second} per degree of freedom, as shown in~\prettyref{subfig:solver_runtimes_homogeneous_times_per_dof}.
Therefore, linear scaling is achieved for~HDG, albeit at relatively high cost.
The transformed system streamlines the operation count of operator and preconditioner so that this quantity is lowered to~\SI{1}{\micro\second}.
This holds for both, continuous and discontinuous case, rendering the choice of discretization one of preference rather than one of speed.
The other variants are somewhat less efficient with larger timings by a factor of to about~$5$.
Still, continuous and discontinuous discretization behave very similar in these cases as well.
The attained hallmark of~\SI{1}{\micro\second} enables high-order solutions to the~\poisson\ equation at costs associated with highly optimized low-order solvers, e.g.\,HPGMG~\cite{gholami_2016_multigrid}, and achieves a main goal of this paper.

\subsection{Robustness against increases in the penalty parameter}

To investigate the robustness against the penalty parameter~$\hdgtauelement$, the tests of the last section were repeated varying the penalty parameter from~$1$, to~$25$, to~$625$.
\begin{figure}[t]
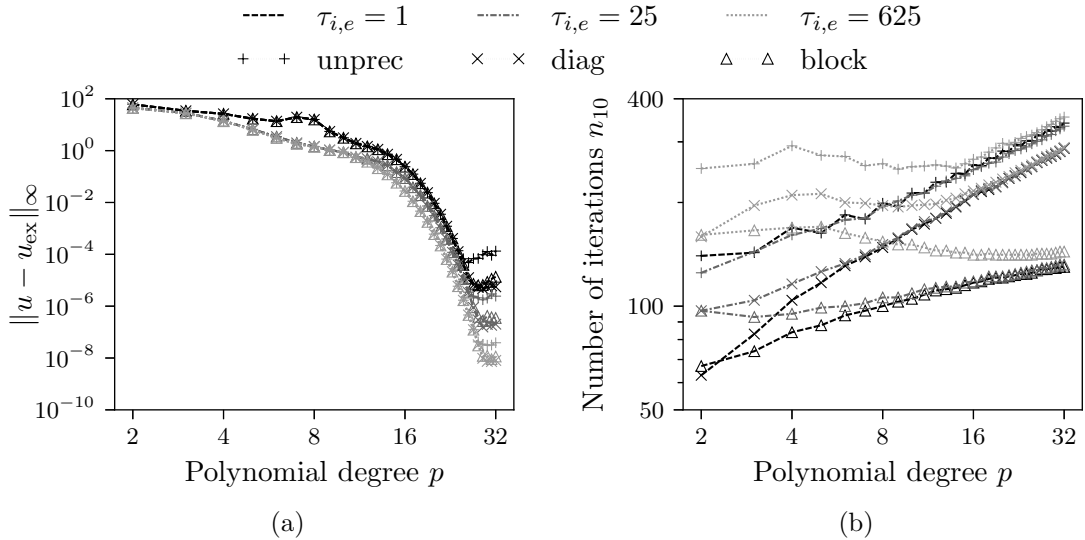

  \hspace*{\fill}
  \includegraphics{hdg_op_plot_condensed_solver_tau_legend.pgf}
  \hspace*{\fill}\\[-1ex]
  \begin{subfigure}{0.49\textwidth}
    \includegraphics{hdg_op_plot_condensed_solver_tau_error_max.pgf}
    \caption{}\label{fig:solver_runtimes_tau_error_max}
  \end{subfigure}
  \begin{subfigure}{0.49\textwidth}
    \includegraphics{hdg_op_plot_condensed_solver_tau_iters.pgf}
    \caption{}\label{fig:solver_runtimes_tau_iters}
  \end{subfigure}
  \caption{Results when varying the penalty parameter~$\hdgtauelement$ for~$\nelement=8^3$ spectral elements of degree~${\poly}$.
    \subref{fig:solver_runtimes_tau_error_max}:~Discretization error.
    \subref{fig:solver_runtimes_tau_iters}:~Number iterations required to lower the residual by ten orders.}%
  \label{fig:solver_runtimes_tau}
\end{figure}%
\prettyref{fig:solver_runtimes_tau} depicts the solution error and the required number of iterations.
An increase of the penalty parameter leads to a reduction of the solution error.
However, the larger penalty parameter takes its toll on the condition of the system.
The influence depends on the polynomial degree:
For~${\poly=2}$, a factor of~$1.5$ occurs when increasing~$\hdgtau$ by a factor of~$25$, whereas for~$\poly=32$, the effect is far lower.
Incidentally, for~$\hdgtauelement = 625$, the condition of the system seems to improve when increasing~$\poly$, again indicating that the system is dominated by the penalty terms at low polynomial degrees.

\subsection{Robustness against increases in the number of elements}

\begin{figure}[t]
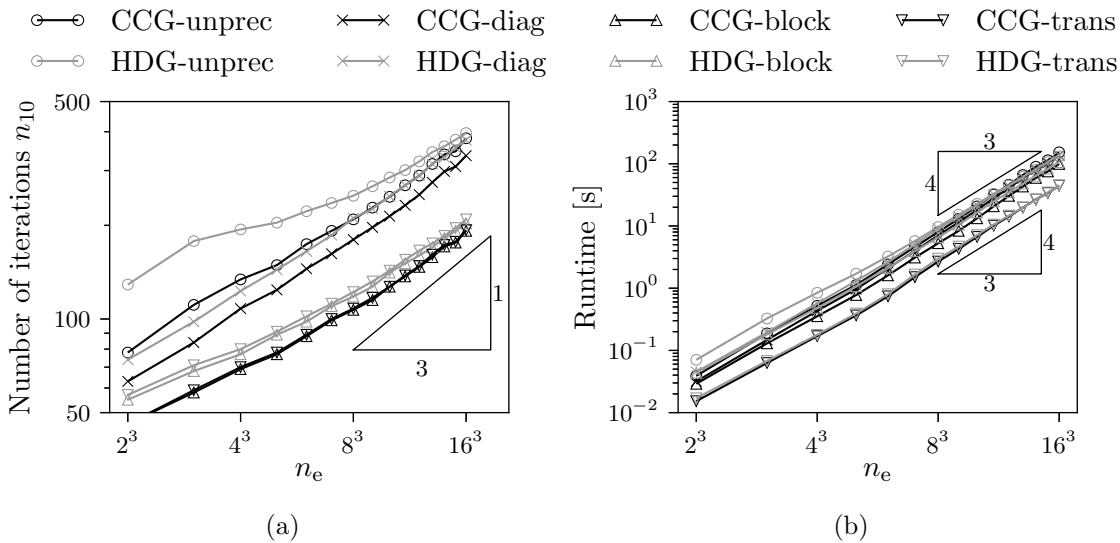

  \hspace*{\fill}
  \includegraphics{hdg_op_plot_condensed_solver_ne_legend.pgf}
  \hspace*{\fill}\\[-1ex]
  \begin{subfigure}{0.49\textwidth}
  \includegraphics{hdg_op_plot_condensed_solver_ne_iters.pgf}
    \caption{}\label{fig:solver_runtimes_ne_iters}
  \end{subfigure}
  \begin{subfigure}{0.49\textwidth}
  \includegraphics{hdg_op_plot_condensed_solver_ne_times.pgf}
    \caption{}\label{fig:solver_runtimes_ne_times}
  \end{subfigure}
  \caption{Results when varying the number of elements~$\nelement$ in a homogeneous mesh using spectral elements of degree~${\poly = 16}$.
    For the~HDG solvers, the penalty parameter was fixed at~${\hdgtauelementi = 25}$.
    \subref{fig:solver_runtimes_ne_iters}:~Number of iterations required to lower the residual by ten orders.
    \subref{fig:solver_runtimes_ne_times}:~Solver runtimes.}%
  \label{fig:solver_runtimes_ne}
\end{figure}
For a constant polynomial degree~${\poly = 16}$, the number of elements is now scaled from~${\nelement = 2^3}$ to~${\nelement = 16^3}$ in this section.
\prettyref{fig:solver_runtimes_ne} depicts the results.
It is evident, that the solvers are not robust against increases in the number of elements, as they exhibit a super-linear increase in runtime with regard to the number of elements.
The slope of the runtime is close to the~$\nelement^{4/3}$ expected for~FEM and finite difference approximations with local preconditioning~\cite{shewchuk_1994_cg}.
To keep the number of iterations constant, global information transport would be required, for instance with multigrid.
In~\cite{stiller_2017a_multigrid}, this was facilitated for~DG by using overlapping~\schwarz-type smoothers.
These are the most promising candidates for extending the present algorithm towards multigrid, and can be factorized to linear complexity for the continuous discretization, as shown in~\cite{huismann_2018_multigrid, huismann_2018_diss}.
However, the derivation of such a factorization for~DG requires enough space for a paper of its own.

\section{Conclusions}\label{sec:conclusions}

This paper considered the residual evaluation for the hybridizable discontinuous~\galerkin\ method for an elliptic equation.
For cuboidal tensor-product elements, the tensor-product decomposition of the operator, in combination with a specific choice of the penalty parameter, and the fast diagonalization technique allowed for a sum factori\-zation.
A linearly scaling evaluation method for the~HDG operator resulted, i.e.\,an operation count proportional to~$\order{\poly^3\nelement} = \order{\ndof}$.
A product factorization lowered the number of operations even further.
The linear scaling was, thereafter, validated with runtime tests.
In these, the operators were as fast as their counterparts for the static condensed continuous case, and even outperformed them for low polynomial orders.
On this basis, linearly scaling Block-\jacobi\ type preconditioners were derived from the operator.
The combination of linearly scaling operator and preconditioner allowed for conjugate gradient solvers with a completely linearly scaling iteration time.
When increasing the polynomial degree, the solvers exhibit only a very slow increase of the number of iterations.
In conjunction with an operator that becomes slightly more efficient at higher polynomial degrees, a near constant iteration time per unknown near~\SI{1}{\micro\second} was achieved.
Not only was this runtime on par with that attained by the continuous solvers, but it also achieves runtimes usually associated with low-order solvers~\cite{gholami_2016_multigrid}.
Compared to these, the solution is by far more accurate so that the proposed methods constitute a leap forward in terms of error reduction compared to the amount of~CPU time invested.

While the solvers proved to scale remarkably well when increasing the polynomial degree, they somewhat lacked robustness against the number of elements.
Therefore, the next step will consist of supplementing them with global coupling.
Multigrid with overlapping~\schwarz\ smoothers is very efficient for both, DG and CG~\cite{haupt_2013_multigrid, stiller_2017_multigrid, stiller_2016_multigrid}.
Moreover, for the static condensed case, the required smoothers were factorized to linear complexity in three dimensions in earlier work~\cite{huismann_2018_multigrid, huismann_2018_diss}, which may now be of help.
However, it is yet unknown whether the approach extends to the discontinuous case.
Iterative substructuring~\cite{sherwin_2001_sem}, the Cascadic multigrid method~\cite{bornemann_1996_cascadic}, or multigrid conjugate gradient methods~\cite{pflaum_2008_multigrid} are further candidates to take into account and were investigated for the continuous methed in~\cite{huismann_2016_cascadic} and~\cite{huismann_2017_substructuring}, respectively.
Work in this direction is under way and will be published in a following paper.

\subsubsection*{Acknowledgements:}
This work is supported in part by the German Research Foundation (DFG) within the Cluster of Excellence ‘Center for Advancing Electronics Dresden’ (cfaed).
The authors to thank their colleagues in the Orchestration path of cfaed for stimulating discussions and ZIH, Dresden, for the computational resources provided.

\bibliographystyle{abbrv}
\bibliography{hdg_op}

\end{document}